\documentclass[journal]{IEEEtran}
\ifCLASSINFOpdf
   \usepackage[pdftex]{graphicx}
\else
\fi
\usepackage{multicol}
\usepackage{multirow}

\usepackage{amsmath}
\usepackage{gensymb}

\usepackage{mathtools}
\usepackage{eurosym}
\DeclareRobustCommand{\officialeuro}{%
  \ifmmode\expandafter\text\fi
  {\fontencoding{U}\fontfamily{eurosym}\selectfont e}}

\usepackage[makeroom]{cancel}
\usepackage[flushleft]{threeparttable}
\ifCLASSOPTIONcompsoc
  \usepackage[caption=false,font=normalsize,labelfont=sf,textfont=sf]{subfig}
\else
  \usepackage[caption=false,font=footnotesize]{subfig}
\fi
\usepackage{url}

\urldef\myurl\url{https://cundall.com/Cundall/fckeditor/editor/images/UserFilesUpload/file/WCIYB/IP-4%20-%20CO2e%20emissions%20from%20biomass%20and%20biofuels.pdf}

\usepackage{amssymb}

\usepackage{nomencl}
\makenomenclature

\usepackage{etoolbox}
\renewcommand\nomgroup[1]{%
  \item[\bfseries
  \ifstrequal{#1}{A}{Index(Sets)}{%
  \ifstrequal{#1}{P}{Parameters}{%
  \ifstrequal{#1}{V}{Variables}{}}}%
]}

\hyphenation{op-tical net-works semi-conduc-tor}

\begin{document}
%
\title{Enforcing Annual Emission Constraints in Short-Term Operation of Local Energy Systems}
%
%
%


%
\author{Dimitri~Pinel\IEEEauthorrefmark{1} and Magnus~Korp\r{a}s\IEEEauthorrefmark{2} \\
Department of Electrical Power Engineering\\
NTNU, Trondheim, Norway\\
E-mail: \IEEEauthorrefmark{1}dimitri.q.a.pinel@ntnu.no, \IEEEauthorrefmark{2} magnus.korpas@ntnu.no\thanks{Address: Elektrobygget, O. S. Bragstads plass 2E, E, 3rd floor, 7034 Trondheim, Norway}}
\maketitle

\begin{abstract}

This paper presents new methods for ensuring that the energy system of a neighborhood that is designed with the objective of being zero emission is actually operated in a way that allows it to reach net zero emissions in its lifetime. This paper highlights the necessity of taking into account realistic operation strategies when designing the energy system of such neighborhoods. It also suggests methods that can be used in the operation of ZENs to ensure carbon neutrality.
An optimization model for designing the energy system of a Zero Emission Neighborhood (ZEN) is first presented and used to produce two designs for a campus in the South of Norway in the case where the amount of PV is limited (PVlim) and when it is not (Base). Several operation approaches are then introduced to compare their operation cost and the $CO_2$ emissions and compensations. These approaches are perfect foresight used as a reference (Ref.), a purely economic model predictive control (E-MPC), an MPC with penalization if deviating from emission targets (EmE-MPC) and a receding horizon MPC where we have a net zero emission constraint over the year (RH-MPC). 
The resulting energy systems are, in the Base case, PV, heat pumps, a gas boiler and heat storage and, in the PVlim case, a smaller amount of PV, a CHP plant, and heat storage. 
In the Base case all operation strategies manage to reach net zero emissions, largely due to the passive compensations obtained from the PV. RH-MPC offers the lowest cost. 
In the PVlim case, the passive effect of the PV is not sufficient to reach net zero emissions and an operation approach specifically taking into account the emissions is necessary. EmE-MPC achieves the lowest emissions but it comes at a much higher cost. We conclude that the best overall strategy is RH-MPC which maintains both the cost and the emission-compensation balance close to the reference case with perfect foresight.

\end{abstract}

\begin{IEEEkeywords}
Operation, Design, Optimization, Distributed Energy Resources, Zero Emission
\end{IEEEkeywords}

\mbox{}
 
\nomenclature[A, 01]{$t (\mathcal{T})$}{Timestep in hour within year, $\in [0,8759]$}
\nomenclature[A, 02]{$\kappa (\mathcal{K})$}{Cluster representative (centroid)}
\nomenclature[A, 03]{$t_\kappa (\mathcal{T_\kappa})$}{Timestep within cluster $\kappa$, $\in [0,23]$}
\nomenclature[A, 04]{$b (\mathcal{B})$}{Building or building type}
\nomenclature[A, 05]{$i (\mathcal{I})$}{Energy technology, $\mathcal{I} = \mathcal{F} \cup \mathcal{E}  \cup \mathcal{HST} \cup \mathcal{EST}; \mathcal{I} = \mathcal{Q} \cup \mathcal{G}$}
\nomenclature[A, 06]{$f (\mathcal{F})$}{Technology consuming fuel (gas, biomass, ...)}
\nomenclature[A, 07]{$e (\mathcal{E})$}{Technology consuming electricity}
\nomenclature[A, 08]{$hst (\mathcal{HST})$}{Heat storage technology}
\nomenclature[A, 09]{$est (\mathcal{EST})$}{Electricity storage technology}
\nomenclature[A, 10]{$q (\mathcal{Q})$}{Technologies producing heat}
\nomenclature[A, 11]{$g (\mathcal{G})$}{Technologies producing electricity}
\nomenclature[A, 12]{$b (\mathcal{B})$}{Building or building type}

\nomenclature[P]{$C_{i,b}^{var,disc}$,$C_{i,b}^{fix,disc}$}{Variable/Fix investment cost of $i$ in $b$ discounted to the beginning of the study including potential re-investments and salvage value [\euro/kWh]/[\euro]}
\nomenclature[P]{$\varepsilon^{tot}_{r,D}$}{discount factor for the duration of the study $D$ with discount rate $r$}
\nomenclature[P]{$C_{i,b}^{maint}$}{Annual maintenance cost of $i$ in $b$  [\euro/kWh]}
\nomenclature[P]{$P_{f}^{fuel}$}{Price of fuel of $g$ [\euro/kWh]}
\nomenclature[P]{$P_{t}^{spot}$}{Spot price of electricity at $t$ [\euro/kWh]}
\nomenclature[P]{$P^{grid}$}{Electricity grid tariff [\euro/kWh]}
\nomenclature[P]{$P^{ret}$}{Retailer tariff on electricity [\euro/kWh]}
\nomenclature[P]{$\eta_{est}$,$\eta_{hst}$}{Efficiency of charge and discharge}
\nomenclature[P]{$\eta_{i}$}{Efficiency of $i$}
\nomenclature[P]{$\eta_{inv}$}{Efficiency of the inverter}
\nomenclature[P]{$\phi_t^{CO_2,el}$}{$CO_2$ factor of electricity at $t$ [$gCO_2/kWh$]}
\nomenclature[P]{$\phi^{CO_2,f}$}{$CO_2$ factor of fuel type $f$ [$gCO_2/kWh$]}
\nomenclature[P]{$\alpha_{CHP}$}{Heat to electricity ratio of the CHP}
\nomenclature[P]{$\alpha_{i}$}{Part load limit as ratio of installed capacity}
\nomenclature[P]{$GC$}{Size of the neighborhood grid connection [kW]}
\nomenclature[P]{$X_{i}^{max}$}{Maximum investment in $i$ [kW]}
\nomenclature[P]{$X_{i}^{min}$}{Minimum investment in $i$ [kW]}
\nomenclature[P]{$E_{b,t}$}{Electric load of $b$ at $t$ [$kWh$]}
\nomenclature[P]{$H^{SH}_{b,t}$,$H^{DHW}_{b,t}$}{Heat (Space heating/Domestic Hot Water) load of $b$ at $t$ [$kWh$]}
\nomenclature[P]{$COP_{hp,b,t}$}{Coefficient of performance of heat pump $hp$}
\nomenclature[P]{$\dot{Q}^{max}_{st}$}{Maximum charge/discharge rate of $est$/$hst$ [kWh/h]}
\nomenclature[P]{$IRR_{t}^{tilt}$}{Total irradiance on a tilted plane [$W/m^2$]}
\nomenclature[P]{$G^{stc}$}{Irradiance in standard test conditions: $1000 W/m^2$}
\nomenclature[P]{$T^{coef}$}{Temperature coefficient}
\nomenclature[P]{$T_{t}$}{Ambient temperature at $t$ [\degree C]}
\nomenclature[P]{$T^{noct}$}{Normal operating cell temperature [\degree C]}
\nomenclature[P]{$T^{stc}$}{Ambient temperature in standard test conditions [\degree C]}
\nomenclature[P]{$\sigma_\kappa$}{Number of occurrences of cluster $\kappa$ in the year}
\nomenclature[P]{$C^{HG}$}{Cost of investing in the heating grid [\euro]}
\nomenclature[P]{$M$}{"Big M", taking a large value}
\nomenclature[P]{$B^{DHW}_{q}$}{Binary parameter stating whether $q$ can produce DHW}
\nomenclature[P]{$Q_{b_1,b_2}^{HGloss}$}{Heat loss in the heating grid in the pipe going from $b_2$ to  $b_1$}
\nomenclature[P]{$\dot{Q}^{MaxPipe}_{b_1,b_2}$}{Maximum heat flow in the heating grid pipe going from $b_2$ to $b_1$ [kWh]}
\nomenclature[P]{$P^{input,max}_{hp,b,t}$}{Maximum power consumption of $hp$ at $t$ based on manufacturer data and output temperature}
\nomenclature[P]{$T^{MPC}$}{Length of the MPC horizon}
\nomenclature[P]{$Em^{0\rightarrow t0}$,$Comp^{0\rightarrow t0}$}{Emissions/Compensations between the start and the current timestep [$gCo_2$]}

\nomenclature[V]{$b^{HG}$}{Binary for the investment in the Heating Grid }
\nomenclature[V]{$b_{i,b}$}{Binary for the investment in $i$ in $b$}
\nomenclature[V]{$x_{i,b}$}{Capacity of $i$ in $b$}
\nomenclature[V]{$f_{f,t,b}$}{Fuel consumed by $f$ in $b$ at $t$ [kWh]}
\nomenclature[V]{$d_{e,t,b}$}{Electricity consumed by $e$ in $b$ at $t$ [kWh]}
\nomenclature[V]{$y_{t}^{imp},y_{t}^{exp}$}{Electricity imported from the grid to the neighborhood/exported at $t$ [kWh]}
\nomenclature[V]{$y_{t,g,b}^{exp}$}{Electricity exported by $g$ to the grid at $t$ [kWh]}
\nomenclature[V]{$g_{t,g,b}^{selfc}$}{Electricity generated by $g$ self consumed in the neighborhood at $t$ [kWh]}
\nomenclature[V]{$g_{t,g,b}^{ch}$}{Electricity generated by $g$ used to charge the batteries at $t$ [kWh]}
\nomenclature[V]{$y_{t,est,b}^{imp}$}{Electricity imported from the grid to $est$ at $t$ [kWh]}
\nomenclature[V]{$y_{t,est,b}^{exp}$}{Electricity exported from the $est$ to the grid at $t$ [kWh]}
\nomenclature[V]{$g_{g,t,b}$}{Electricity generated by $g$ at $t$ [kWh]}
\nomenclature[V]{$q_{q,t,b}$}{Heat generated by $q$ in $b$ at $t$ [kWh]}
\nomenclature[V]{$y^{dch}_{t,est,b}$}{Electricity discharged from $est$ to the neighborhood at $t$ [kWh]}
\nomenclature[V]{$y^{ch}_{t,est,b}$}{Electricity charged from on-site production to $est$ at $t$ [kWh]}
\nomenclature[V]{$q^{ch}_{t,st,b},q^{dch}_{t,st,b}$}{Energy charged/discharged from the neighborhood to the storage at $t$ [kWh]}
\nomenclature[V]{$v^{stor}_{t,st,b}$}{Level of the storage $st$ in building $b$ at $t$  [kWh]}
\nomenclature[V]{$g^{curt}_{t,b}$}{Solar energy production curtailed [kWh]}
\nomenclature[V]{$g^{dump}_{g,t,b}$}{Electricity generated but dumped by $g$ at $t$ [kWh]}
\nomenclature[V]{$q^{dump}_{t,b}$}{Heat dumped at $t$ in $b$ [kWh]}
\nomenclature[V]{$q^{HG transfer}_{b_1,b_2,t}$}{Heat transferred via the heating grid from $b_1$ to $b_2$ at $t$ [kWh]}
\nomenclature[V]{$q^{HG used}_{b,t}$}{Heat taken from the heating grid by $b$ at $t$ [kWh]}
\nomenclature[V]{$o_{i,t,b}$}{Binary controlling if $i$ in $b$ is on or off at $t$}
\nomenclature[V]{$\overline{x_{i,b,t}}$}{Maximum production from $i$ [kWh]}
\nomenclature[V]{$c^{Em}$,$c^{Comp}$}{Penalization cost for deviating from the emission/compensation targets [\euro]}

\printnomenclature[2cm]

%
\IEEEpeerreviewmaketitle

\section{Introduction}

The control of energy systems usually only takes into account the near future. Indeed, the plan of operation of a system can only be as good as the forecast fed to it, and such forecasts' precision quickly drops with the length of the forecast horizon. In many cases, this does not create problems. Let us take the example of a battery energy storage system in a house. For a control algorithm based on Model Predictive Control (MPC), i.e. an optimization with a rolling horizon where only the first time step is implemented, we can model the way the battery should be operated based on the spot price of electricity and load forecasts in order to minimize electricity costs.
However, if the model is also intended to capture trends affecting the operation of the battery on a longer time scale than your horizon, the way to incorporate these trends to the model is not straightforward. Such a trend could for instance be a grid tariff design with peak power pricing or calendar aging of the battery in the operation. This paper investigates possible ways to include long-term trends in the specific case of the control of Zero Emission Neighborhoods' energy system.

Zero Emission Neighborhoods (ZENs) are neighborhoods that aim to have no net emissions of $CO_2$ in their lifetime. In order to design the energy system of such neighborhoods, a tool called ZENIT has been developed. It uses a mixed integer programming (MIP) optimization to minimize the cost of investing in and operating the energy system of a ZEN. However, the way the system is operated in the investment optimization is important in order to reach this net zero emission criteria during actual operation. How can we make sure that the system is operated in a similar way and that the long-term goal of zero net emissions is captured by the control scheme of the neighborhood? This paper explores this problem by making propositions of ways to handle this issue and evaluating them.

As will be shown in section \ref{sec2}, the existing literature fails to address the problematic of incorporating long-term targets used in the design process into the short-term operation of systems. This paper contributes to the existing literature on neighborhood energy system planning and MPC operation by:
\begin{itemize}
    \item highlighting the gap between designs of the energy system of neighborhoods and their actual operation
    \item introducing novel methods for handling short term-operation while incorporating long-term goals, specifically net zero emission requirements
    \item analyzing the performance of each method from an operational point of view and from the point of view of the goal of net zero emission
\end{itemize}
It deals with identifying the pros and cons of different approaches for a practical implementation of ZENs.

In section \ref{sec2}, the existing literature regarding the operation of the energy system of neighborhoods is presented and the research gaps highlighted. In section \ref{stat_ana}, the zero emission objective is presented and the features of the input timeseries used in the models for different years are compared in order to select an appropriate reference year for the analysis. The implications of the zero emission goal for the operation are discussed. In section \ref{inv_mod} the investment model is briefly introduced and the resulting neighborhoods that will be used in the rest of the study are presented. Following this, in section \ref{op_mod}, the different model alternatives are presented, and their results analysed in section \ref{op_res}.

\section{State of the Art and Contribution} \label{sec2}

The investment in the energy system of more sustainable neighborhoods and buildings or even ZENs is an area where several models have been proposed.
\cite{MEHLERI2012}, for example, considers this problem together with the heating grid layout problem and includes emissions in the objective function of the MILP via a carbon tax.
Recently, \cite{GABRIELLI2018} and \cite{FLEISCHHACKER2019} suggested multi-objective approaches, where the objectives of cost and emission minimization are opposed to obtain a Pareto front of solutions. \cite{GABRIELLI2018} models seasonal storages and addresses the difficulties arising from using clustering in this context. \cite{FLEISCHHACKER2019} uses two open source models (urbs and rivus) together with different levels of spatial and temporal aggregation through clustering; one is for the design of the energy system of the buildings and the other is for the design of the electrical and heating grid.

MPC emerged in the 1960s and gained traction in the oil and chemical industries \cite{MORARI1999667}. The principle of MPC is to control a process through a model-based on-line optimization strategy. At each timestep, an optimization with a finite horizon is solved in order to obtain the plan of operation of the next period. Only the first timestep is actually implemented. The optimization is then run again for the next period, moving one timestep forward, taking into account the actual realization of the actions resulting from the previous run and the updated forecast.
It is now used in various applications. Controlling the energy systems of buildings or neighborhoods is such an example \cite{OLDEWURTEL201215,Ma12, SIROKY20113079, oldewurtel10, PRIVARA2011564, MOROSAN20101445, halvgaard12, maasoumy12, Celik18, MA201292}. These works focus on indoor climate control and/or energy management (load shifting, peak shaving). Those schemes typically use a thermal model of the building they operate in order to model the inside temperature and comfort level.

Other ways of operating an energy system exist, for example rule-based operations or systems operated manually by users only, but they do not offer the same level of comfort or efficiency \cite{LEE2016}.

From the perspective of planning the energy system of a neighborhood, the actual control of the system cannot be included in detail in the models due to complexity reasons and simplifications are thus necessary. The investment optimization will use a deterministic or stochastic optimal operation and the actual system can then be operated in various ways that are likely to be sub-optimal in the long run, affecting only the short-term operation cost. Both the investment and the operation are done in a purely economic way.

From the point of view of a planner of a ZEN, this is different. There is a constraint to have zero net emissions of $CO_2$ during its lifetime. This requires a specific operation of the system and an accounting of the emitted and compensated $CO_2$. To the best of the authors knowledge, no paper deals with the inclusion of integral constraint, or long-term-goals, in modelling the on-line operation of the  energy system of a neighborhood. The authors are not aware of similar problems in different contexts, with the exception of hydropower scheduling and battery degradation. In \cite{BODAL2019} for example, the long-term hydro-reservoir level curves resulting from long-term planning models of the reservoirs are used to constrain the production of the hydro plant in a rolling horizon framework. This approach is standard in the operation of hydropower, with the use and linkage of long-, medium- and short-term models. In \cite{ortega14}, the degradation of the battery over its lifetime is accounted for in its daily operation.
In the context of emission accounting, a white paper from Soteica \cite{soteica2009}, the implementation of an integral constraint on emission of different pollutant is discussed including $CO_2$ emissions. The context is different than in this paper. They consider the operation of a refinery where the emission constraint coming from the allowance given in the EU ETS (Emission Trading System) needs to be taken into account in the daily operation. They use different approaches depending on the remaining emissions from the allowance. Four approaches are discussed briefly in \cite{soteica2009}: 
\begin{itemize}
    \item assigning a fixed cost to the emissions
    \item no cost until the quota is achieved
    \item assigning a cost to the emission based on the projections for the emissions level in the rest of the year
    \item an emission limit based on the same projections
\end{itemize}

The context of a ZEN is different, since there is no allowance, but compensations are also accounted for and a balance between the emissions and compensations is required. A different approach is necessary in the ZEN context even if some ideas from \cite{soteica2009} can be applied.

This review of exiting literature highlights the need for better ways to deals with long-term goals in the short term planning. Our contribution is to introduce several new methods to deal with this problem. To this end we investigate these approaches in an operation model of a ZEN energy system and use the results of the investment model as a reference. In the following section we will present the zero emission objective and analyze the features of different historical years to select the most representative reference year used in the investment model.

\section{Selection of the Reference Year} \label{stat_ana}

\subsection{Zero Emission Objective}

For our neighborhood to be a ZEN, we need to meet the zero emission requirement. This means that the neighborhood should have net zero emissions at the end of its lifetime. What should be included in the emissions of the neighborhood varies depending on the ambition of the stakeholders. It can simply be the emissions from the operation but can also include the embedded emissions of the material, the emissions from the construction and the deconstruction of the neighborhood. In order to reach zero emissions, the emissions need to be compensated. In this study we only focus on the emissions coming from the operation phase of the neighborhood's lifetime.
The approach in ZENIT is to consider that the export of electricity from on-site renewable generation sources to the grid reduces the amount of electricity produced nationally, contributing to emissions with a higher $CO_2$ factor. The emissions that were avoided thanks to the export from the neighborhood are accounted as compensation in the zero emission balance. 

In ZENIT, the optimization model uses one representative year for the lifetime of the neighborhood. In order to give good insights into the necessary investments, the reference year should have average electricity price and temperature levels. The temperatures should also represent minimum temperatures correctly because this will have an effect on the maximum heat demand. In order to ensure a good representation of the compensations from PV it should also have average solar conditions.

When operating a neighborhood that was designed to become a ZEN, the question is if you should try to have a zero emission balance every year. Indeed, what was the case for the reference year is not necessary for specific years, a year with lower than average solar irradiances could, for instance, be compensated by a year with higher than average irradiances and the emission even out over the lifetime. The different methods proposed in section  \ref{op_mod} try to impact the operation by considering the reference year's emissions and for that reason they do not take this possibility of year-to-year compensation into account. This constitutes a shortcoming of the presented methods.

\subsection{Statistical analysis of the Inputs over the Years}

To get meaningful results, it is important to consider wisely the year to be used in the optimization model. The choice of year can impact the results significantly \cite{PFENNINGER2017} \cite{jafari2020}. With a limited availability of data, it is therefore important to consider the available years carefully.

In order to determine the appropriate reference year, but also in order to know the features of the input data for different years and be able to analyze the results in the rest of the paper we present the boxplot, duration curve and density curve of various input timeseries. The inputs selected are the outside temperature, the solar irradiance, the spot price of electricity and the $CO_2$ factor of electricity. The loads of the buildings are not included because we assume a strong correlation to the outside temperature. The years included are 2015, 2016, 2017 and 2018 because they are the years for which we have the timeseries of $CO_2$ factors of electricity.

\begin{figure}
    \centering
    \includegraphics[width=0.48\textwidth]{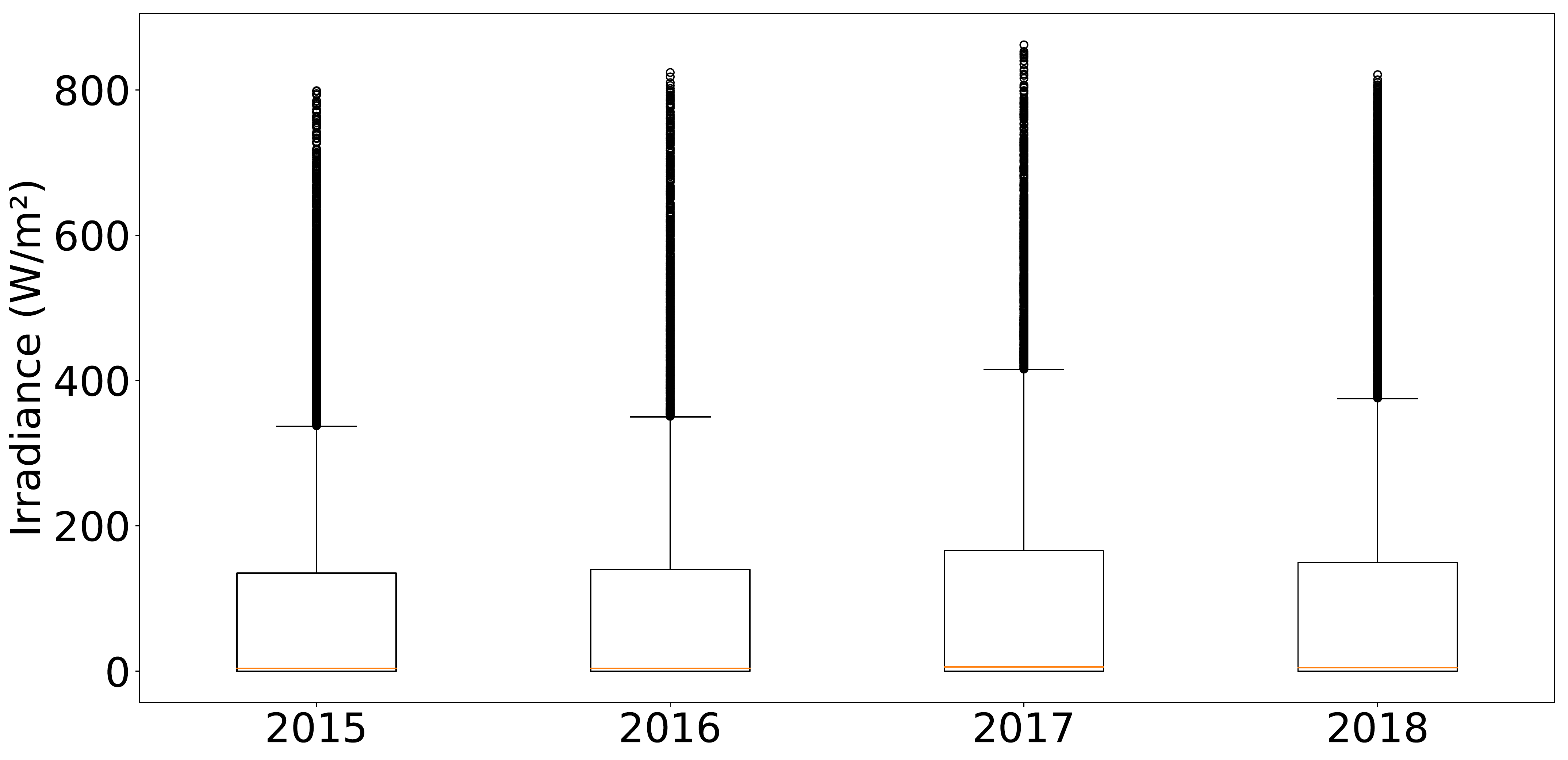}
    \caption{Boxplot of Solar Irradiance for the Studied Years}
    \label{fig:box_irr}
\end{figure}

\begin{figure}
    \centering
    \includegraphics[width=0.48\textwidth]{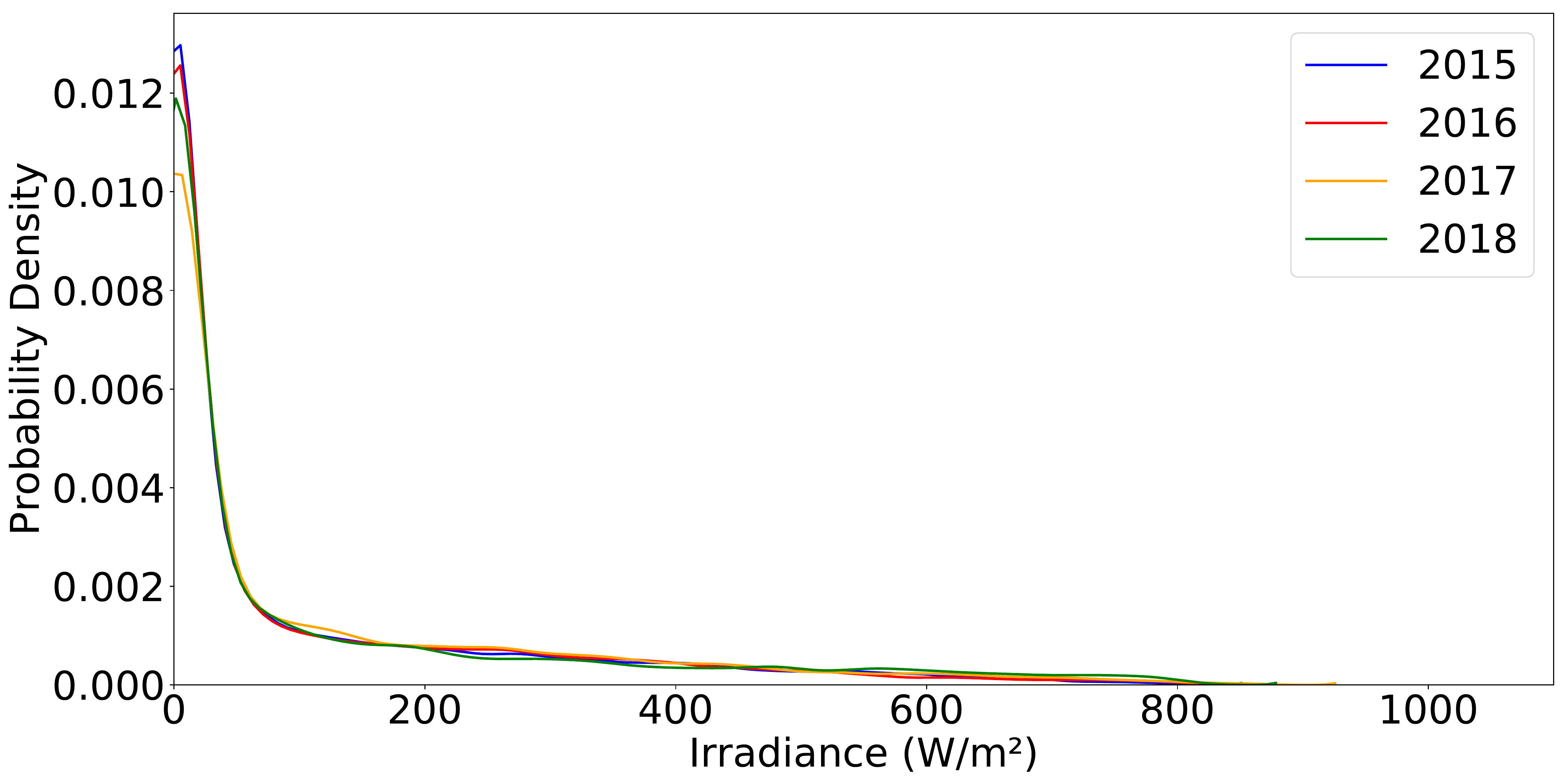}
    \caption{Density Curves of Solar Irradiance for the Studied Years}
    \label{fig:den_irr}
\end{figure}

\begin{figure}
    \centering
    \includegraphics[width=0.48\textwidth]{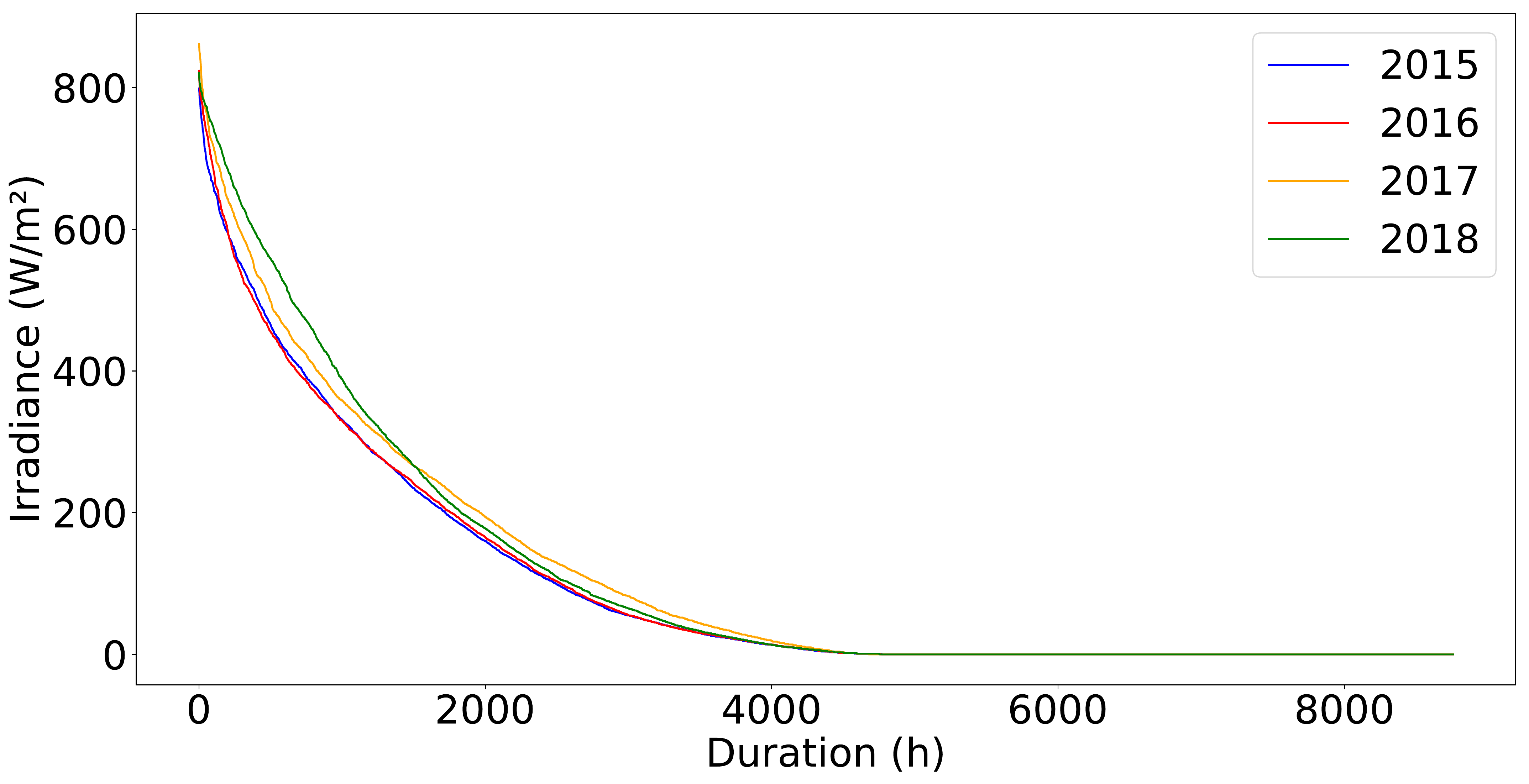}
    \caption{Duration Curves of Solar Irradiance for the Studied Years}
    \label{fig:dur_irr}
\end{figure}

The solar irradiance is quite similar for the different years with minor variations due to weather conditions. The years 2018 and 2017 have the highest total irradiance. From the density curves Fig. \ref{fig:den_irr}, we can see that there is roughly the same probability for the irradiance to be between 0 and 100 than above $100W/m^2$. The boxplot Fig. \ref{fig:box_irr} confirms the distribution, with a median close to 0, a third quartile around 130 and numerous outliers.

\begin{figure}
    \centering
    \includegraphics[width=0.48\textwidth]{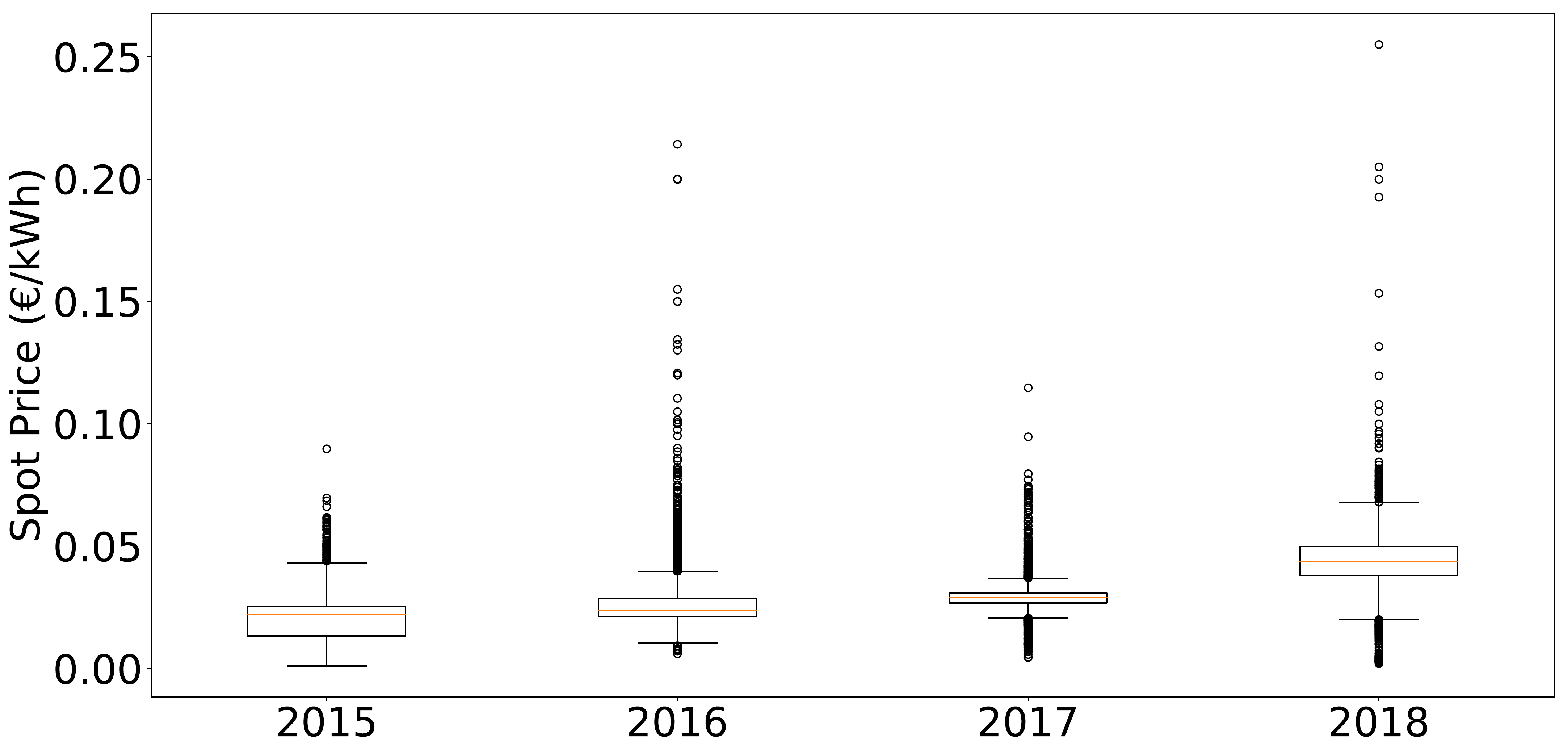}
    \caption{Boxplot of Spot Prices for the Studied Years}
    \label{fig:box_spot}
\end{figure}

\begin{figure}
    \centering
    \includegraphics[width=0.48\textwidth]{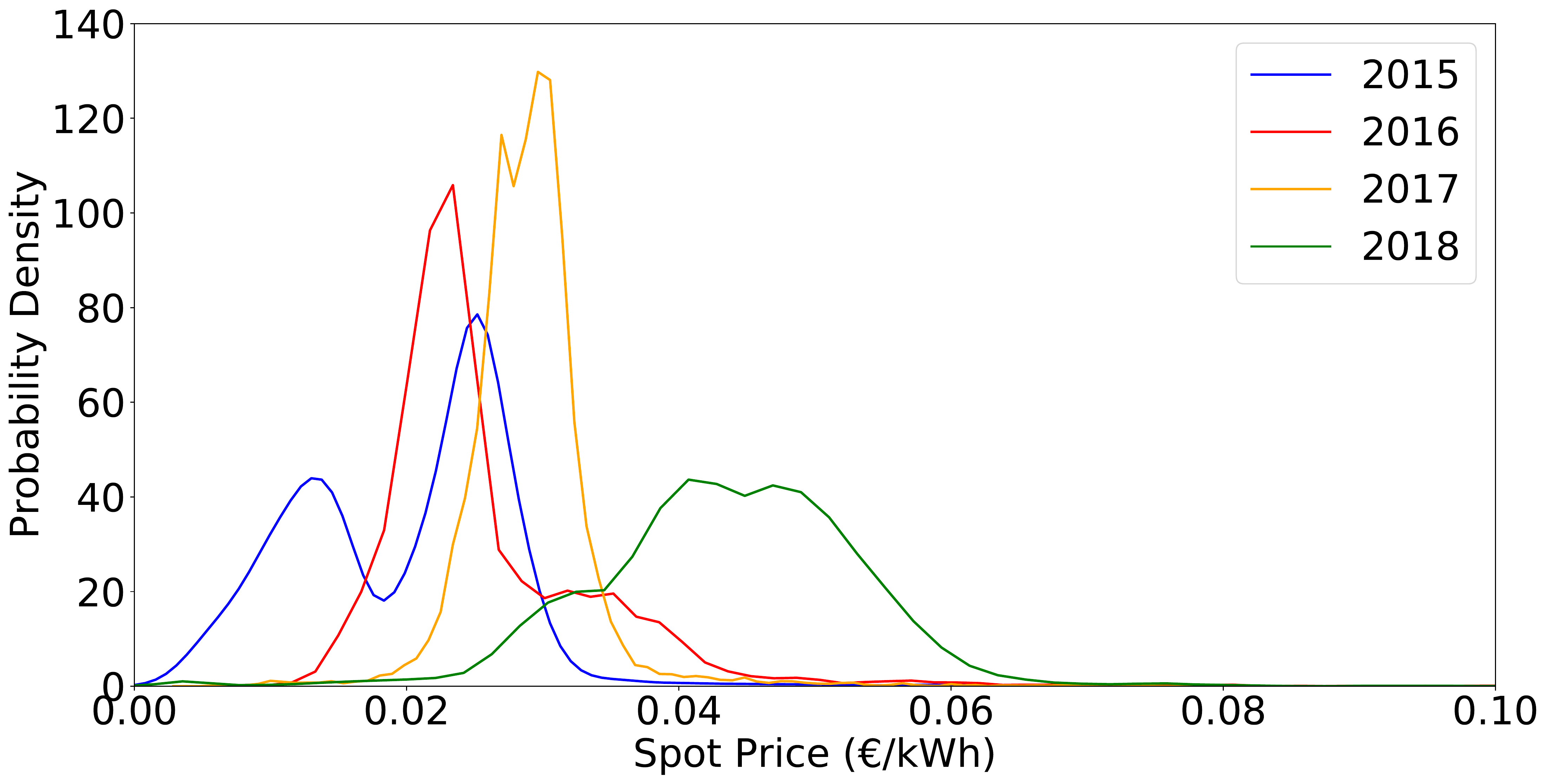}
    \caption{Density Curves of Spot Prices for the Studied Years}
    \label{fig:den_spot}
\end{figure}

\begin{figure}
    \centering
    \includegraphics[width=0.48\textwidth]{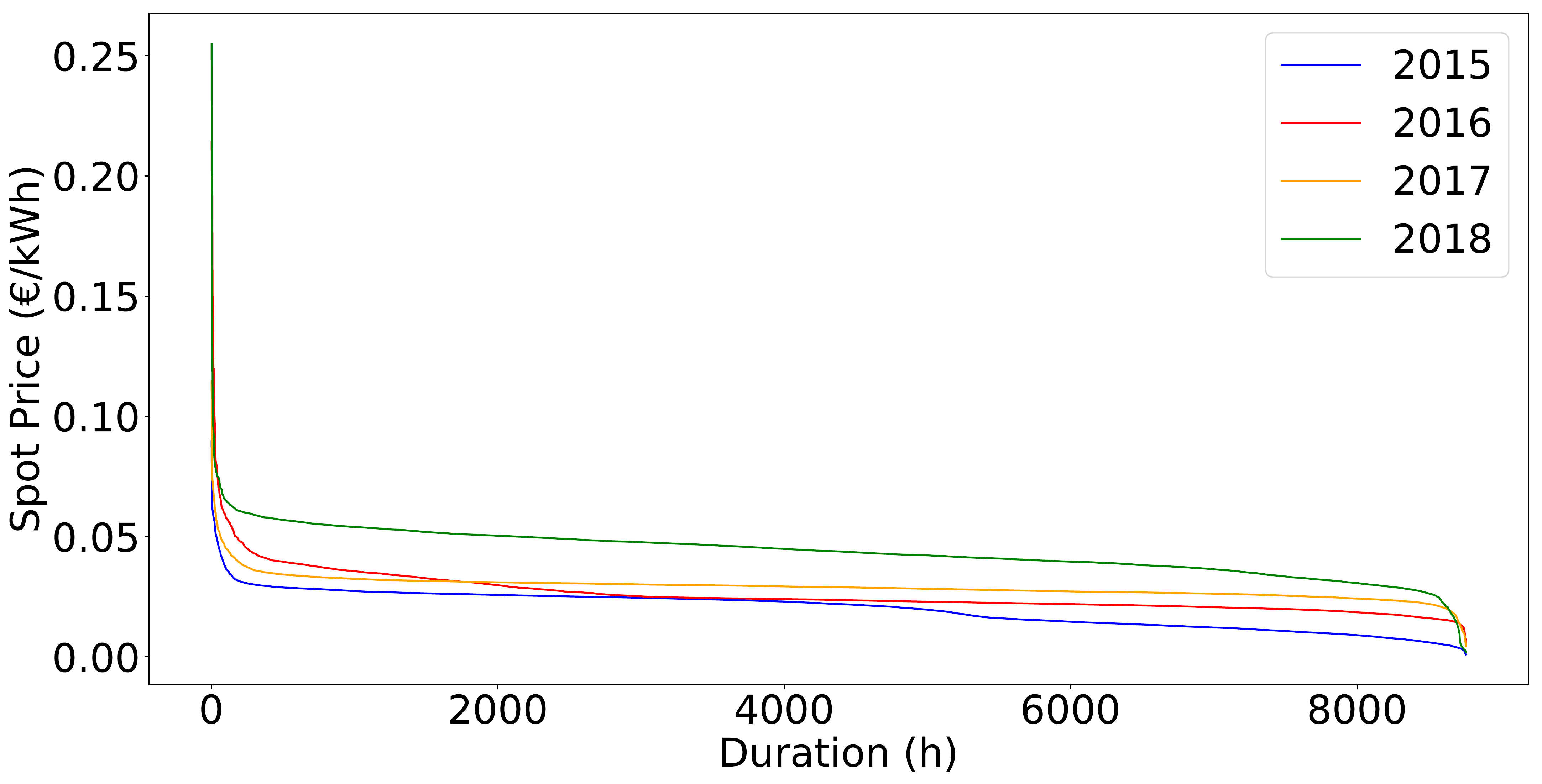}
    \caption{Duration Curves of Spot Prices for the Studied Years}
    \label{fig:dur_spot}
\end{figure}

From Fig. \ref{fig:box_spot}, we can see that the values of spot prices are not noticeably spread, the bands between quartile 1 and 3 are narrow. However, there are some outliers, mainly reflecting peaks in prices but also dips for the case of 2017 and 2018. The median values also vary significantly. It is also important to note the difference in highest peak prices in 2016 and 2018 compared to 2017 and 2015. The distribution of the prices shown in Fig. \ref{fig:den_spot} are quite different. They are all relatively wide with the exception of 2017, but the shape and the means are quite different. The year 2018, for instance, is more even while the rest have a peak, denoting the concentration of the prices around that value. In the case of 2015, there are two peaks denoting two price levels where most of the data lie. Those observations are confirmed by the duration curve Fig. \ref{fig:dur_spot}.

\begin{figure}
    \centering
    \includegraphics[width=0.48\textwidth]{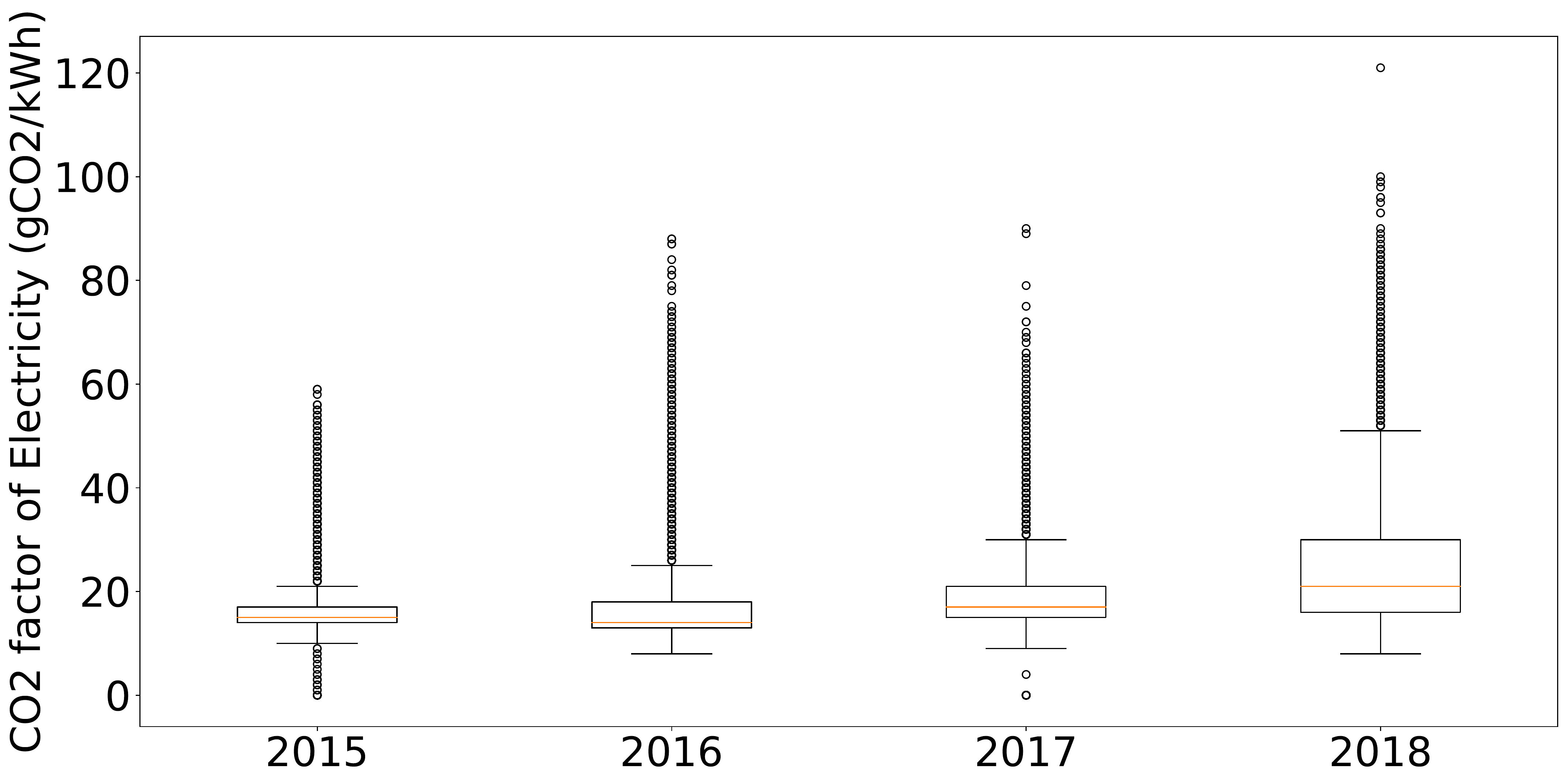}
    \caption{Boxplot of $CO_2$ factors of electricity for the Studied Years}
    \label{fig:box_co2}
\end{figure}

\begin{figure}
    \centering
    \includegraphics[width=0.48\textwidth]{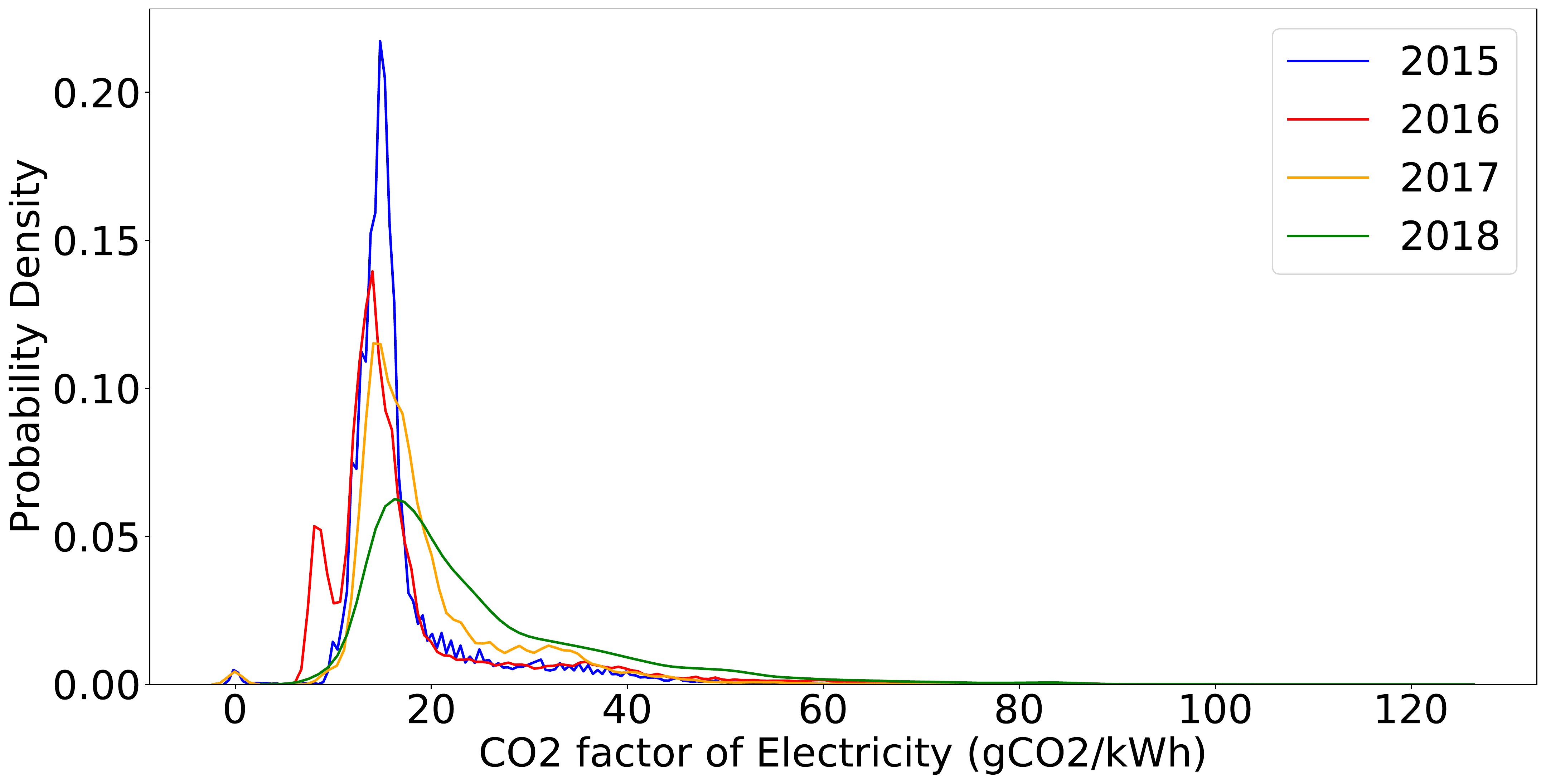}
    \caption{Density Curves of $CO_2$ factors of electricity for the Studied Years}
    \label{fig:den_co2}
\end{figure}

\begin{figure}
    \centering
    \includegraphics[width=0.48\textwidth]{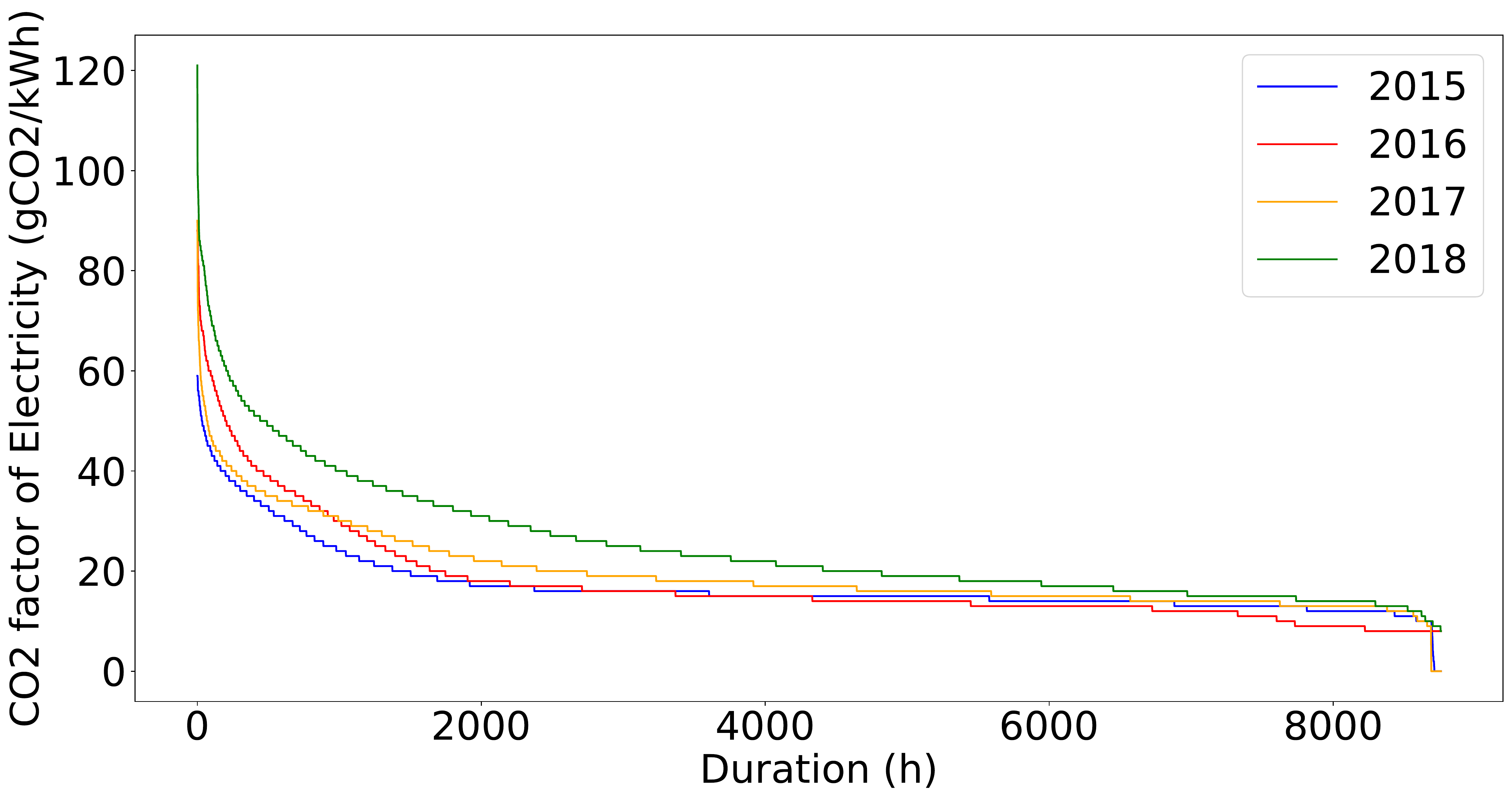}
    \caption{Duration Curves of $CO_2$ factors of electricity for the Studied Years}
    \label{fig:dur_co2}
\end{figure}

The $CO_2$ factors for electricity also show 2018 as quite different from other years, with higher median (Fig. \ref{fig:box_co2}) and wider distribution of values (Fig. \ref{fig:den_co2}). The other years are more similar with a median of around 17 $gCO_2/kWh$. The year 2016 offers a somewhat middle-ground representation of the peak levels of the $CO_2$ factors even if the base levels are slightly lower than for other years (Fig. \ref{fig:dur_co2}).

\begin{figure}
    \centering
    \includegraphics[width=0.48\textwidth]{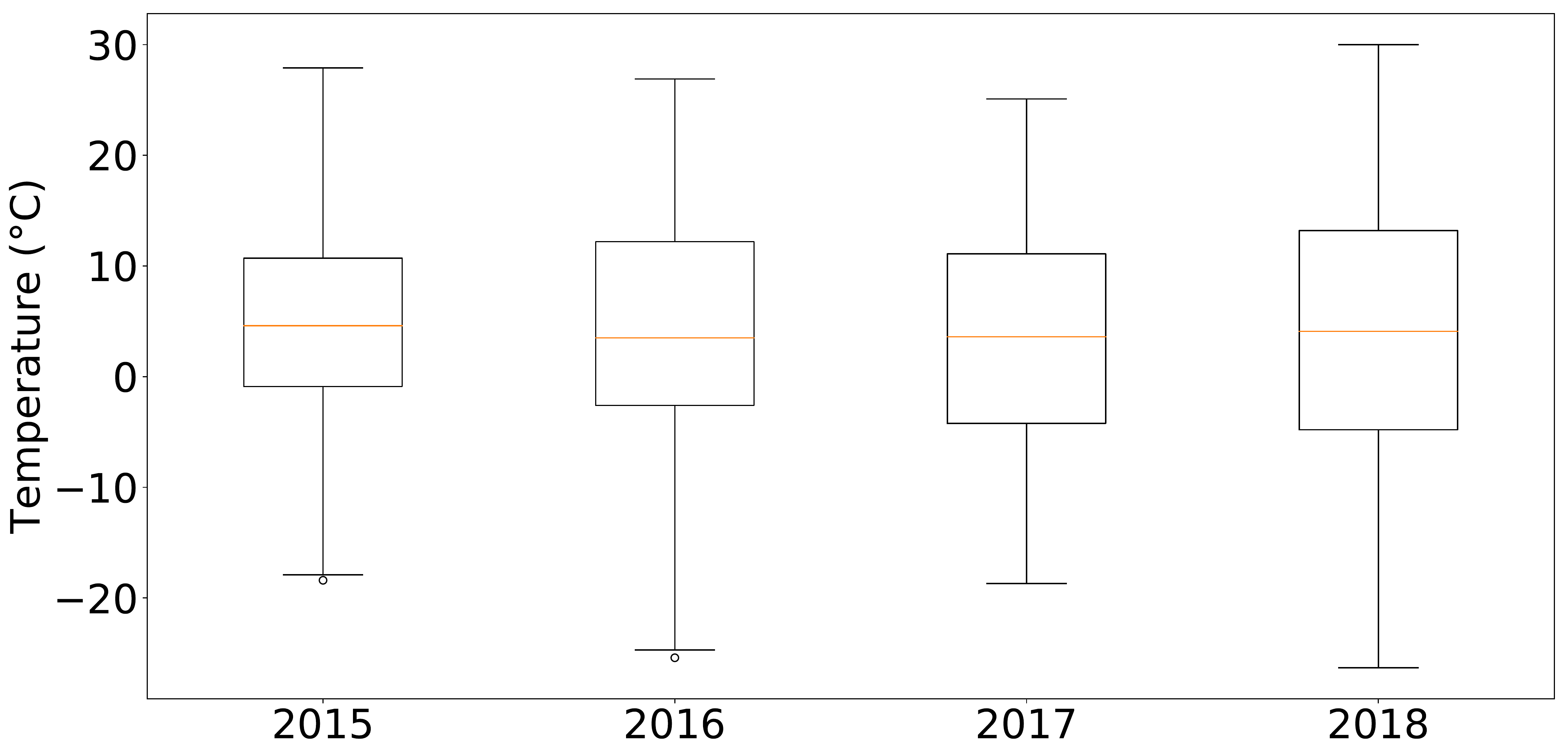}
    \caption{Boxplot of Temperature for the Studied Years}
    \label{fig:box_temp}
\end{figure}

\begin{figure}
    \centering
    \includegraphics[width=0.48\textwidth]{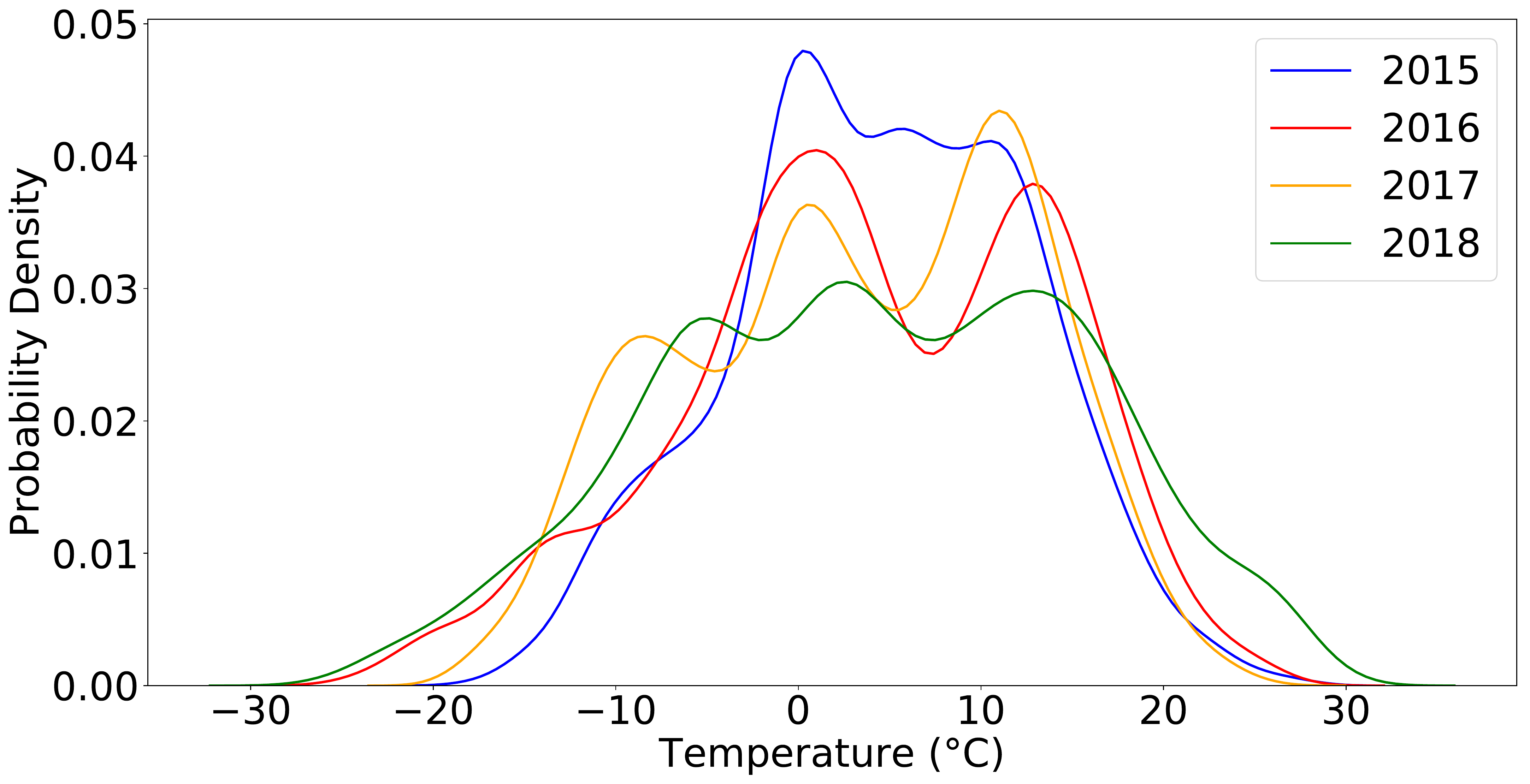}
    \caption{Density Curves of Temperature for the Studied Years}
    \label{fig:den_temp}
\end{figure}

\begin{figure}
    \centering
    \includegraphics[width=0.48\textwidth]{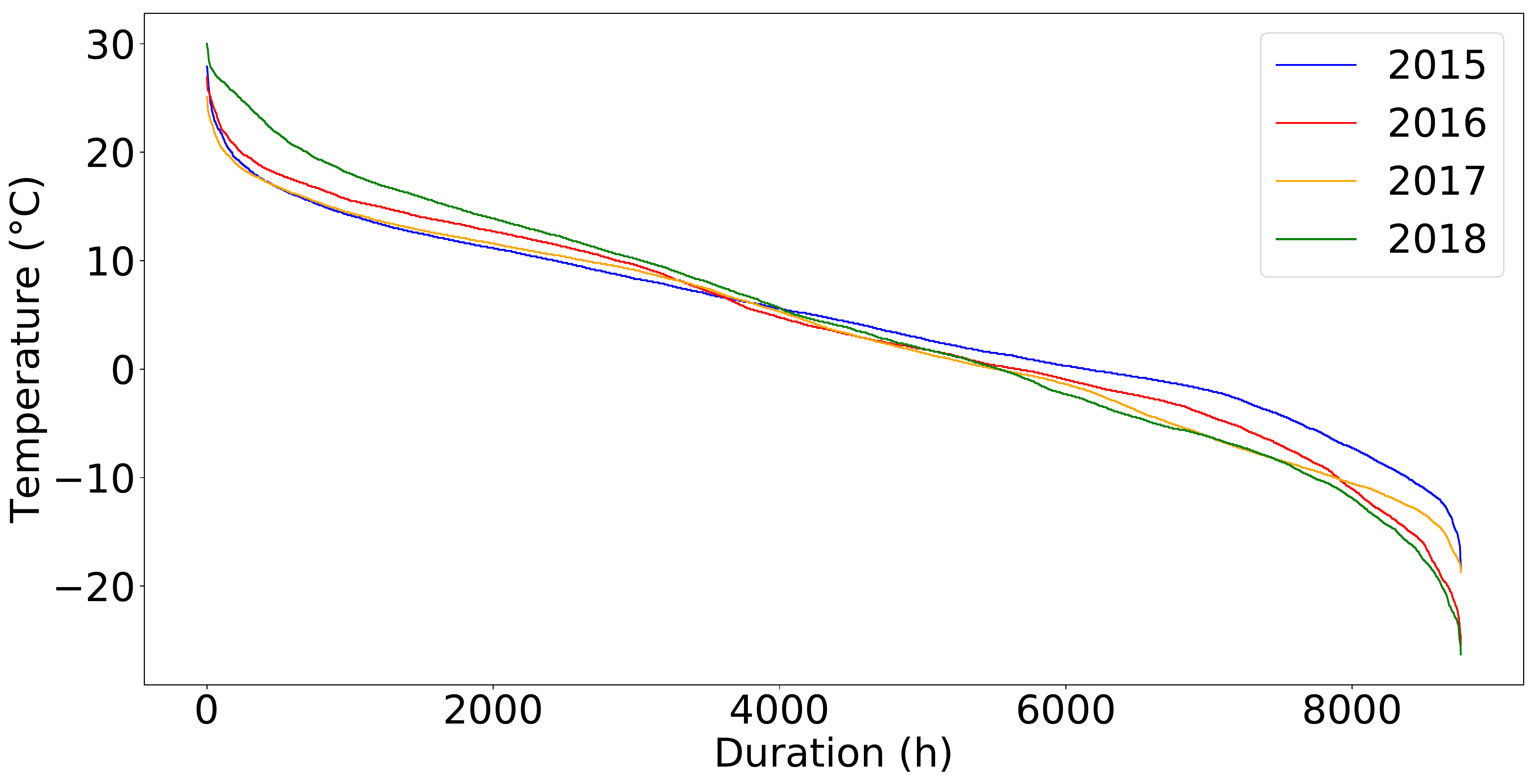}
    \caption{Duration Curves of Temperature for the Studied Years}
    \label{fig:dur_temp}
\end{figure}

The median of the temperature lies around $5\degree C$ for all years, as seen in Fig. \ref{fig:box_temp}. There is a bigger spread of values than for the other timeseries and almost no outliers. The distribution of the different years Fig. \ref{fig:den_temp} is quite similar even if their shape varies.

Overall it seems that 2016 is a good candidate to be used as a reference year for the investment optimization from our sample of years. It has average temperatures while still having high and low extremes (Fig. \ref{fig:dur_temp}). It also has a somewhat average representation of the solar irradiance and of the spot price. The representation of the $CO_2$ factors is also average for the ``peaks" but slightly lower in the base level. We chose this year to make the investment optimization for these reasons and keep the three other years to use them in the comparison of the operation strategies.

\section{Investment Model ZENIT and its Results} \label{inv_mod}

\subsection{Presentation of ZENIT} \label{zenit}
In this section, the investment model called ZENIT (Zero Emission Neighborhood Investment Tool) and the setup of the study are presented before introducing the resulting systems that will be operated in the following sections.

ZENIT uses optimization to find the cost-optimal energy system for a neighborhood to be zero emission. It uses one representative year instead of the whole lifetime for computational reasons. This description is an extract from \cite{dimitri19}.
The objective function is:

\textit{Minimize}:
\begin{multline}
    b^{HG}\cdot C^{HG} + \sum_{b}\sum_{i}  \Big( (C_{i,b}^{var,disc} + \frac{C_{i,b}^{maint}}{\varepsilon^{tot}_{r,D}})\cdot x_{i,b} +\\ C_{i,b}^{fix,disc} \cdot b_{i,b} \Big) 
     + \sum_{t_\kappa}\frac{\sigma_\kappa}{\varepsilon^{tot}_{r,D}}\Big(\sum_b\sum_{f}f_{f,t,b} \cdot P_{f}^{fuel}\\
     + (P_{t}^{spot}+P^{grid}+P^{ret})\cdot (y_{t}^{imp} \\
     +\sum_b\sum_{est}y_{t,est,b}^{imp})-P_{t}^{spot}\cdot y_{t}^{exp}\Big)
\end{multline}

It considers the fix and variable investment cost of the different technologies ($C_{i,b}^{var,disc}$, $C_{i,b}^{fix,disc}$) and the heating grid ($C^{HG}$), as well as operation- and maintenance-related costs ($C_{i,b}^{maint}$). A binary variable controls the investment in the heating grid ($b^{HG}$). The subscripts used in the equations are $b$ for the buildings, $i$ for the technologies, $t$ for the timesteps, $f$ for fuels and $est$ for batteries. $\varepsilon$ are the discount factors with interest rate $r$ for the duration of the study $D$. $x_{i,b}$ is the capacity of the technologies and $b_{i,b}$ the binary related to whether it is invested in or not. $\sigma_\kappa$ is the number of occurrences of cluster $\kappa$ in the full year and $t_\kappa$ is the timestep in the cluster. $P$ are the prices of fuel, electricity on the spot market, grid tariff or retailer tariff. $f$ is the consumption of fuel and $y$ are the imports or exports of electricity.

In order to fulfill the zero emission requirement presented in section \ref{stat_ana}, the following constraint, called the Zero Emission Balance is used:

\begin{multline}\label{zebal}
    \phi^{CO_2,el}_{t}\sum_{t_\kappa}\sigma_\kappa\Big(y_{t}^{imp}+\sum_b\sum_{est}y_{t,est,b}^{imp}\Big) \\
    + \sum_{t_\kappa}\sigma_\kappa\sum_b\sum_{f} \phi^{CO_2,f} \cdot f_{f,t,b} \leq  \sum_{t_\kappa} \phi^{CO_2,el}_{t}\cdot\sigma_\kappa \\
    \bigg( \sum_b\sum_{est}\eta_{est}\cdot y_{t,est,b}^{exp} 
    +\sum_b\sum_{g} y_{t,g,b}^{exp}\bigg)
\end{multline}

The $CO_2$ factors are represented by $\phi^{CO_2,el}_{t}$ for electricity and $\phi^{CO_2,f}$ for other fuels. $\eta_{est}$ is the charging efficiency of the battery.

Other equations include load balances for electricity (\ref{elec_bal}), domestic hot water (DHW) (\ref{DHW_bal}) and space heating (SH) (\ref{SH_bal}). They require the production and import to be equal to the consumption and exports for all timesteps. $\forall t$:

\begin{subequations}
\begin{multline} \label{elec_bal}
    y^{imp}_{t}+\sum_b\Big(\sum_{est}y^{dch}_{t,est,b} \cdot \eta_{est} +\sum_{g} g_{g,t,b}^{selfc}\Big) \\= \sum_b\Big(\sum_e d_{e,t,b} + E_{b,t} \Big)
\end{multline}
$\forall t,b$:
\begin{multline} \label{DHW_bal}
    \sum_{q} q^{DHW}_{q,t,b}+ \sum_{hst}(\eta_{hst} \cdot q^{DHWdch}_{t,hst,b}-q^{DHWch}_{t,hst,b})\\+q^{HGusedDHW}_{t,b} = H^{DHW}_{b,t}  + q^{dump}_{t,b}
\end{multline}
\begin{multline} \label{SH_bal}
    \sum_{q} q^{SH}_{q,t,b}+ \sum_{hst}(\eta_{hst} \cdot q^{SHdch}_{t,hst,b}-q^{SHch}_{t,hst,b})\\+q^{HGusedSH}_{t,b} =  H^{SH}_{b,t} 
\end{multline}
\end{subequations}

The optimization model can choose to invest in a heating grid (\ref{hgpt}), giving access to other technologies. We assume that those technologies are located in a central production plant that feeds the heating grid. The operation of the heating grid is then constrained by the following equations: $\forall t$
\begin{subequations}
\begin{multline} \label{energy bal pp}
    \sum_{q} q_{q,t,'PP'}+ \sum_{hst}(\eta_{hst} \cdot q^{dch}_{t,hst,'PP'}-q^{ch}_{t,hst,'PP'}) \\= \sum_{b\setminus'PP'}q^{HGtrans}_{t,'PP',b} + q^{dump}_{t,'PP'}
\end{multline}
$\forall b,b',t$
\begin{equation} \label{flow lim}
    q_{t,b',b}^{HGtrans} \leq \dot{Q}^{MaxPipe}_{b',b}
\end{equation}
$\forall b,t$
\begin{equation}\label{out limit}
    \sum_{b'} q_{t,b,b'}^{HGtrans} \leq \sum_{b''} \Big( q_{t,b'',b}^{HGtrans} - Q_{b'',b}^{HGloss} \Big)
\end{equation}
\begin{equation} \label{hg Sh DHW}
    q_{t,b}^{HGused}=q_{t,b}^{HGusedSH}+q_{t,b}^{HGusedDHW}
\end{equation}
\begin{equation}\label{heat from HG in b}
    q_{t,b}^{HGused}= \sum_{b''} \Big( q_{t,b'',b}^{HGtrans} - Q_{b'',b}^{HGloss} \Big) - \sum_{b'} q_{t,b,b'}^{HGtrans}
\end{equation}
$\forall i$ 
\begin{equation}\label{hgpt}
     x_{i,'Production Plant'} \leq X_{i}^{max} \cdot b^{HG}
\end{equation}
\end{subequations}

The energy balance at the central production plant (PP in the equations) is modelled with \ref{energy bal pp}, the flow limit in the pipes by \ref{flow lim}, the distinction between the heat from the heating grid used for SH or DHW by \ref{hg Sh DHW}, and the heat used in the specific building by \ref{heat from HG in b}. Equation \ref{out limit} sets the maximum for what goes out of the building to what came in, i.e. heat produced in the building cannot be fed to the heating grid.

The connection to the national electric grid limits the exports and imports: $\forall t$
\begin{equation}
    y_{t}^{imp}+\sum_b\sum_{est}y_{t,est,b}^{imp}+\sum_b\sum_{g}y_{t,g,b}^{exp} \leq GC
\end{equation}

For most technologies, the production of heat or electricity is linked to the fuel consumption using the efficiency of the technology.
\begin{subequations}
$\forall \gamma \in \mathcal{F} \cap \mathcal{Q},t,b$:
\begin{equation}
    f_{\gamma,t,b}=\frac{q_{\gamma,t,b}}{\eta_\gamma}
\end{equation}
$\forall \gamma \in \mathcal{E} \cap \mathcal{Q},t,b$:
\begin{equation}
    d_{\gamma,t,b}=\frac{q_{\gamma,t,b}}{\eta_\gamma}
\end{equation}
\end{subequations}

For CHPs the electricity produced is the ratio of the heat produced and the heat to power ratio $\alpha_{CHP}$: $\forall t,'CHP',b$:
\begin{equation}
    g_{CHP,t,b}= \frac{q_{CHP,t,b}}{\alpha_{CHP}}
\end{equation}

The heat produced can be used for DHW or for SH (\ref{shdhw}) but some technologies can only provide SH (such as electric radiators or wood stove). Equation \ref{onlysh} translates this constraint.
$\forall q,t,b$:
\begin{equation}\label{shdhw}
    q_{q,t,b} = q^{DHW}_{q,t,b} + q^{SH}_{q,t,b}
\end{equation}
\begin{equation}\label{onlysh}
    q^{DHW}_{q,t,b} <= M \cdot B^{DHW}_q
\end{equation}

The production from PV and solar thermal collectors depends on the irradiance on a tilted surface $IRR_{t}^{tilt}$ and their efficiency. The efficiency for the solar panel $\eta^{PV}_t$ is defined based on \cite{hellman14} and accounts for the cell temperature $T_c$ and inverter losses.
\begin{subequations}
\begin{equation}
    g_{PV,t}+g^{curt}_{t} = \eta_{PV,t} \cdot x_{PV} \cdot IRR_{t}^{tilt}
\end{equation}
\begin{equation}
    q_{ST,t}= \eta_{ST} \cdot x_{ST} \cdot IRR_{t}^{tilt}
\end{equation}
\begin{equation}
    \eta_{PV,t}=\frac{ \eta^{inv}}{G^{stc}} \cdot \big(1-T^{coef}\cdot(T^c-T^{stc})\big)
\end{equation}
\begin{equation}
    T^c=T_t+(T^{noct}-20)\cdot \frac{IRR_t^{tilt}}{800}
\end{equation}
\end{subequations}
For the heat pumps in the buildings, the production and electrical consumption are defined as follows:
\begin{subequations}
\begin{equation}\label{eq:copsh}
    d_{hp,b,t}^{SH}=\frac{q_{hp,b,t}^{SH}}{COP_{hp,b,t}^{SH}}
\end{equation}
\begin{equation}\label{eq:copdhw}
    d_{hp,b,t}^{DHW}=\frac{q_{hp,b,t}^{DHW}}{COP_{hp,b,t}^{DHW}}
\end{equation}
\begin{equation}\label{eq:HP_lim}
    \frac{d_{hp,b,t}^{DHW}}{P^{input,max,DHW}_{hp,b,t}} + \frac{d_{hp,b,t}^{SH}}{P^{input,max,SH}_{hp,b,t}} \leq x_{hp,b}
\end{equation}
\end{subequations}

Equations \ref{eq:copsh} and \ref{eq:copdhw} link the heat produced to the COP and the electrical consumption of the heat pump. The COPs are different for SH and DHW due to different temperature set points. They also depend on the outside temperature and they are calculated before the optimization.
Equation \ref{eq:HP_lim} regulates how the heat pump can be used for both SH and DHW and enforces that the capacity invested is not exceeded. $P^{input,max}$ represents the maximum power input to the heat pump at the timestep based on the temperature set point and for a 1kW unit. $d_{hp,b,t}^{SH}$ and $d_{hp,b,t}^{SH}$ represent the electric consumption of the heat pump for SH and DHW while $q_{hp,b,t}^{DHW}$ and $q_{hp,b,t}^{DHW}$ are the heat production.

Another binary variable is used for part load limitations. This binary concerns the operation and is defined for every timestep for each relevant technology, which can lead to a large number of binary variables. No minimum up- or downtime is used. $\forall i\setminus HP,t,b$:
\begin{subequations}
\begin{equation}
    \overline{x_{i,b,t}} \leq X_{i,b}^{max} \cdot o_{i,t,b}
\end{equation}
\begin{equation}
    \overline{x_{i,b,t}} \leq x_{i,b}
\end{equation}
\begin{equation}
    \overline{x_{i,b,t}} \geq x_{i,b} - X_{i,b}^{max} \cdot (1-o_{i,t,b})
\end{equation}
\begin{equation}
   q_{i,b,t} \leq \overline{x_{i,b,t}} 
\end{equation}
\begin{equation}
   q_{i,b,t} \geq \alpha_{i,b} \cdot \overline{x_{i,b,t}} 
\end{equation}
\end{subequations}
The size of the investment in each technology type is bounded from below to represent the larger scale of some technologies (\ref{mininv}) and from above (\ref{maxinv}) to limit the size of the research space. $\forall i,b$:
\begin{subequations}
\begin{equation}\label{mininv}
    x_{i,b} \leq X_{i,b}^{max} \cdot b_{i,b}
\end{equation}
\begin{equation} \label{maxinv}
    x_{i,b} \geq X_{i,b}^{min} \cdot b_{i,b}
\end{equation}
\end{subequations}

Technologies producing electricity can feed this electricity to the neighborhood directly, store it in batteries, export it or dump it.$\forall t,g,b$:
\begin{equation}
    g_{g,t,b}=y_{t,g,b}^{exp}+g_{g,t,b}^{selfc}+g_{t,g,b}^{ch}+g^{dump}_{t,g,b}
\end{equation}

To distribute the production to the batteries, we have $\forall t,b$:
\begin{equation}
    \sum_g g_{t,g,b}^{ch} = \sum_{est} y^{ch}_{t,est,b}
\end{equation}

The storage operation, be it heat or electrical storage, is modeled as follows:
$\forall \kappa, t_\kappa \in [1,23],st,b$
\begin{equation}
    v^{stor}_{\kappa,t_\kappa,st,b}=v^{stor}_{\kappa,t_\kappa-1,st,b}+\eta^{stor}_{st,b}\cdot q_{\kappa,t_\kappa,st,b}^{ch} -q_{\kappa,t_\kappa,st,b}^{dch}
\end{equation}

$\forall \kappa, t_\kappa \in [0,23],st,b$
\begin{equation}
     v^{stor}_{\kappa,t_\kappa,st,b} \leq x_{st,b}
\end{equation}

\begin{multicols}{2}\noindent
    \begin{equation}
     q^{ch}_{\kappa,t_\kappa,st,b} \leq \dot{Q}_{st}^{max}
    \end{equation}\noindent
    \begin{equation}
     q^{dch}_{\kappa,t_\kappa,st,b} \leq \dot{Q}_{st}^{max}
    \end{equation}
\end{multicols}
$\forall st,b,\kappa$
\begin{equation}\label{eq:endval}
    v^{stor}_{\kappa,0,st,b}=v^{stor}_{\kappa,23,st,b}
\end{equation}

The state of charge of the storage $st$ (either heat or electric storage) is represented by $v^{stor}$ while $q^{ch}$ and $q^{dch}$ are the energy charged and discharged. The maximum charge and discharge rate is ${Q}_{st}^{max}$. This model only allows for the use of representative days and daily storage operation. Details of the process of clustering and choice of an appropriate number of clusters can be found in \cite{dimitri19}.

We perform two investment runs. In the first one the roof area constrains the amount of solar technologies that can be installed. In the second one we assume that there is available area in the proximity that can be used to install solar panels and we do not take the roof area into account.

The model is implemented on a test case based on a small neighborhood, a campus at Evenstad in Norway, where three building types represent the different buildings there. We use the same implementation as in \cite{dimitri19}. More information on the implementation of the studied case can be found there. 

The $CO_2$ factors for electricity are obtained by tracing back the origin of the electricity using the methodology presented in \cite{clauss19}. The data used in this methodology primarily comes from the ENTSO-E transparency platform. The earliest complete data on the platform start in 2015, which explains our choice of years.

The investment options details and sources are presented in Annex \ref{annex:data}.

\subsection{Results from ZENIT} \label{sec_inv_res}

The results from the investment runs are presented in this subsection. 
In the rest of the paper we will refer to ``Base" and ``PVlim" for, respectively, the case where PV is not constrained by PV area and the case where it is. The central plant represents the location where the neighborhood scale technologies are and we refer to the existing buildings from the Campus Evenstad as B1, B2 and B3. B1 represents student apartments at the passive standard, B2 conventional offices and B3 offices at the passive standard.

\begin{figure}[h]
    \centering
    \includegraphics[width=0.48\textwidth]{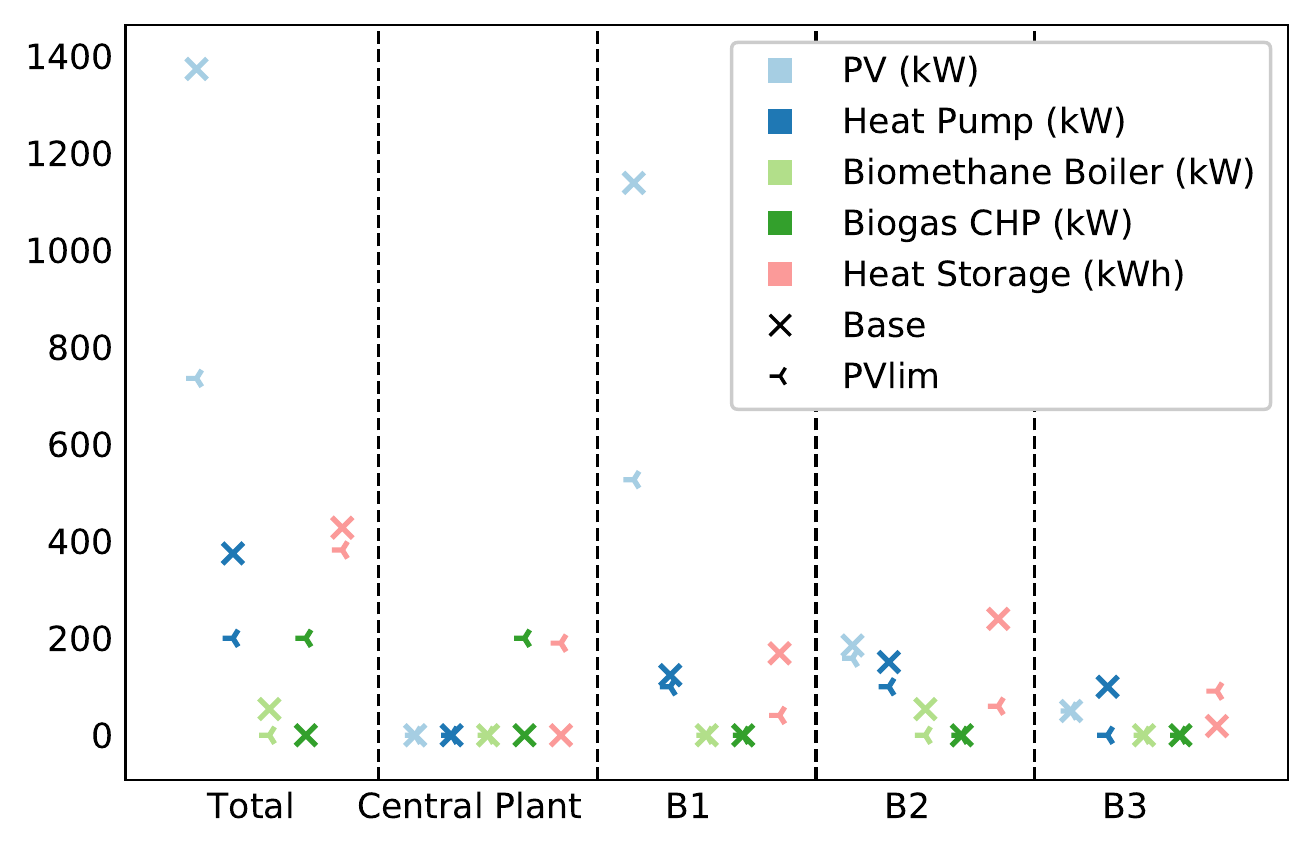}
    \caption{Results of the Investment Runs in the Base and PVlim case}
    \label{inv_res}
\end{figure}

The investments resulting from the runs for the Base and PVlim cases are shown in Fig. \ref{inv_res}.
In the Base case, a combination of a large amount of PV, air-source heat pumps are used together with a biomethane boiler in B2.
In the PVlim case, the amount of installed PV is around two times lower than in the Base case. The limitation on PV also induces an investment in the heating grid and a biogas engine at the neighborhood level. This partially replaces heat pumps in particular in B3 and completely replaces the biomethane boiler in B2.

The total emissions for one year are respectively 11.69 and 5.25 $tonCO_2$ for the Base and PVlim case. The compensations are strictly equal to the emissions.

\begin{table}[h]
    \centering
    \caption{Discounted Investment and Discounted Operation Costs for the Base and PVlim Cases in k\euro; the Sum Represents the Optimal Objective Value}
    \begin{tabular}{l r r r}
    \hline
         & Disc. Investment Cost & Disc. Operation Cost\\
         \hline
        Base & 1 351.1 & 993.7\\
        PVlim & 1 077.3 & 1 706.7 \\
        \hline
    \end{tabular}
    \label{tab:inv_costs}
\end{table}

The discounted investment and operation costs are shown in Table \ref{tab:inv_costs}. The lifetime used for the neighborhood is 60 years and the rate of return is 4\%.

\section{Operation Models} \label{op_mod}

In this section we present the different models used to assess the operation of the neighborhood.

\subsection{Reference Model}

The reference model operates the neighborhood with perfect foresight. It is able to operate the neighborhood in a perfect way and is thus used as a reference value for the other methods. This is however not a method that can be used in practice due to the increasing errors of forecasts with time.

This also represents the way the system would have been operated by the investment optimization. Indeed, we use the same formulation for the optimization with the exception that the investment part is removed. The objective function becomes:

\textit{Minimize}:
\begin{multline}
    \sum_{t=0}^{8759}\Big(\sum_b\sum_{f}f_{f,t,b} \cdot P_{f}^{fuel}
     + (P_{t}^{spot}+P^{grid}+P^{ret})\\\cdot (y_{t}^{imp} 
     +\sum_b\sum_{est}y_{t,est,b}^{imp})-P_{t}^{spot}\cdot y_{t}^{exp}\Big)
\end{multline}

\subsection{Economic MPC (E-MPC)}

The model that we call economic MPC or E-MPC uses the same constraints as the reference model but uses a rolling horizon of 24 hours to operate the system. There is no perfect foresight anymore and the operation thus cannot anticipate future conditions of prices or temperatures for example. One optimization is run for each timestep and only the first timestep is implemented. Since there is no actual operation of a system there is no problem regarding the difference between the plan for a timestep and the actual realization for this timestep, which means we assume that the operation plan decided by the optimization is perfectly realized.
The objective function becomes: 

$\forall t0 \in [0..8759]:$
\textit{Minimize}:
\begin{multline}
    \sum_{t=t0}^{t0+T_{MPC}}\Big(\sum_b\sum_{f}f_{f,t,b} \cdot P_{f}^{fuel}
     + (P_{t}^{spot}+P^{grid}+P^{ret})\\\cdot (y_{t}^{imp} 
     +\sum_b\sum_{est}y_{t,est,b}^{imp})-P_{t}^{spot}\cdot y_{t}^{exp}\Big)
\end{multline}
With $T_{MPC}$ the length of the horizon, which is 24 hours in our case. The constraints stay the same as in the previous models, except that they are defined over the horizon only. The operation of the storages links the different horizons through the storage level at $t0$. 

\subsection{Emission Constrained Economic MPC (EmE-MPC)}

The emission constrained MPC (EmE-MPC) uses the same formulation as for the E-MPC but adds a penalization cost for deviating from emission and compensation targets. The targets are calculated based on the results from the investment runs. One emission target and one compensation target is calculated for each horizon. The penalization is added to the objective function, which becomes:
$\forall t0 \in [0..8759]:$
\textit{Minimize}:
\begin{multline}
    \sum_{t=t0}^{t0+24}\Big(\sum_b\sum_{f}f_{f,t,b} \cdot P_{f}^{fuel}
     + (P_{t}^{spot}+P^{grid}+P^{ret})\cdot (y_{t}^{imp} \\
     +\sum_b\sum_{est}y_{t,est,b}^{imp})-P_{t}^{spot}\cdot y_{t}^{exp}\Big) + c^{Em} + c^{Comp}
\end{multline}

\begin{figure}
    \centering
    \includegraphics[width=0.48\textwidth]{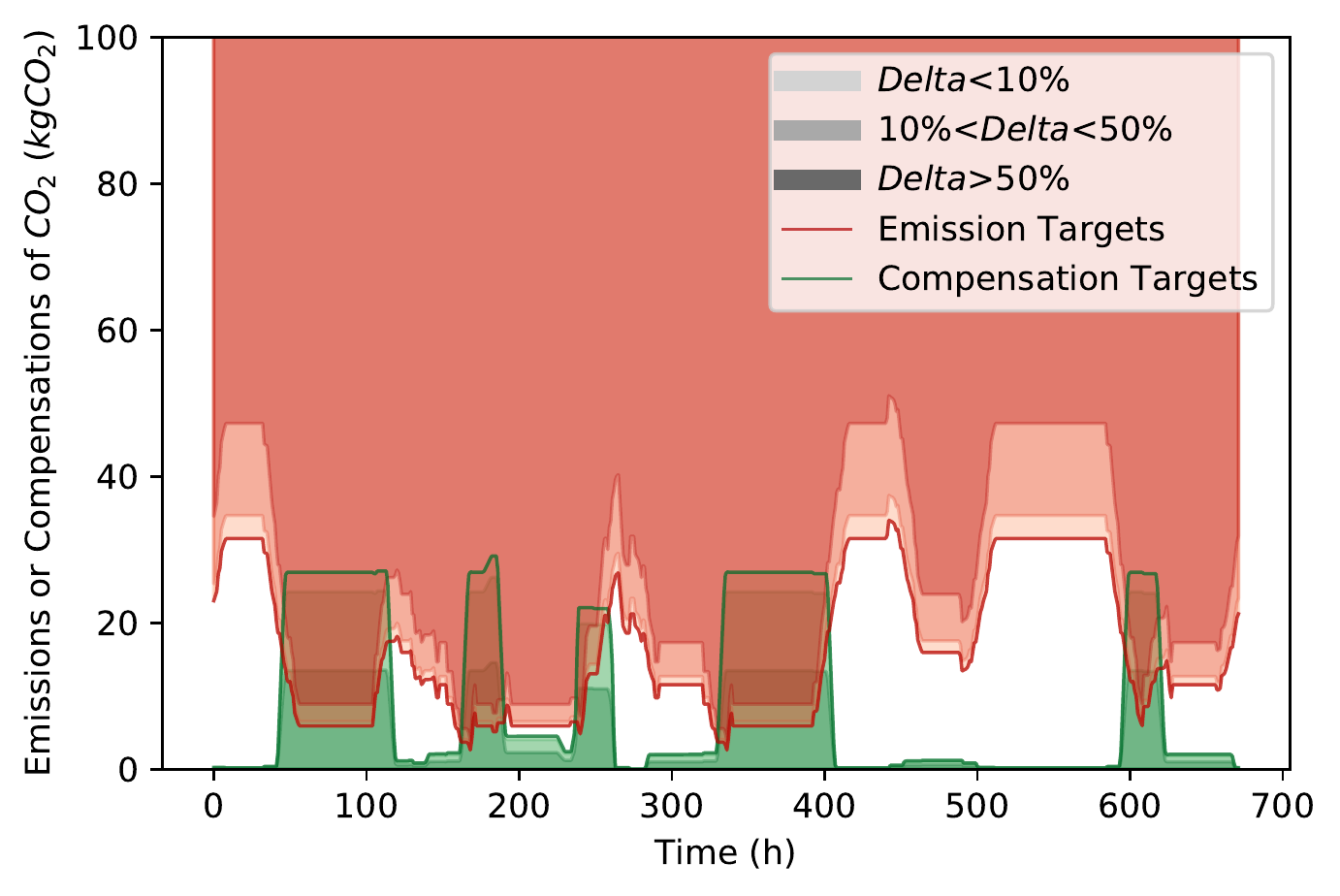}
    \caption{Representation of the Different Emission and Compensation Targets Ranges (Delta: Difference Between Actual and Target) used in Equations \ref{eq:em} and \ref{eq:comp} for Each Horizon for One Winter Month for PVlim Case}
    \label{fig:em_comp_target}
\end{figure}

The penalization is calculated in the following way:
\begin{equation} \label{eq:em}
    c^{Em}=\delta_1\cdot(Em^{1.1})+\delta_2\cdot(Em^{1.5})+\delta_3\cdot(Em^{sup})
\end{equation}
\begin{multline} \label{eq:comp}
    c^{Comp}=\delta_3\cdot(Comp^{0})\cdot b^{0}+\delta_2\cdot(Comp^{0.5})\cdot b^{0.5}\\+\delta_1\cdot(Comp^{0.9})\cdot b^{0.9}
\end{multline}
Where $Em^{1.1}$ are the emissions up to $10\%$ above the emission target, $Em^{1.5}$ the emissions between $10$ and $50\%$ above the emission target and $Em^{sup}$ the emissions above the latter. 
For the compensation, the calculation is different and has discontinuities. $Comp^{0}$, $Comp^{0.5}$ and $Comp^{0.9}$ are the difference between the compensation target and the actual compensation when this compensation is respectively between 0 and 50$\%$, 50 and 90$\%$, and 90 and 100$\%$ of the target value.
Figure \ref{fig:em_comp_target} represents the emissions and compensations targets and ranges for each horizon in one winter month in the PVlim case. The areas above the red line and below the green line are respectively the ranges of penalized emissions and compensations. The white area either below or above represents values of emissions or compensations that represent less emissions or more compensations and as such do not get penalized.

Only one of the three components at most is active in the equation because of the binaries and the following equation ($b^{sup}$ represents the case of compensations higher than the target):
\begin{equation}
    b^{0}+b^{0.5}+b^{0.9}+b^{sup}=1
\end{equation}
The emissions and compensations are calculated with the same formulas as respectively the left-hand side and the right-hand side of the zero emission balance, equation \ref{zebal}.

The values of $\delta_1$, $\delta_2$ and $\delta_3$ were set after multiple tries to respectively 0.03, 3 and 300$\euro/gCO_2$.

\subsection{Receding Horizon MPC (RH-MPC)}

In the receding horizon MPC (RH-MPC) we use a complete year so that we are able to re-introduce the zero emission balance over the year. To maintain similar foresight conditions as in the previous models, we use the timeseries values of the next horizon only from the actual year to operate and we use the reference year values for the rest of the year. From $t0$ to $t0+t_{MPC}$ the corresponding data in the current year are used and for $t0+t_{MPC}$ to 8759 we use the reference year data.

The objective function becomes:
$\forall t0 \in [0..8759]:$
\begin{multline}
    \sum_{t=t0}^{8759}\Big(\sum_b\sum_{f}f_{f,t,b} \cdot P_{f}^{fuel}
     + (P_{t}^{spot}+P^{grid}+P^{ret})\cdot (y_{t}^{imp} \\
     +\sum_b\sum_{est}y_{t,est,b}^{grid\_imp})-P_{t}^{spot}\cdot y_{t}^{exp}\Big) 
\end{multline}

The emission balance constraint is reintroduced in the following form:
\begin{multline}
    Em^{0\rightarrow t0}+\phi^{CO_2,el}_{t}\sum_{t=t0}^{8759}\Big(y_{t}^{imp}+\sum_b\sum_{est}y_{t,est,b}^{imp}\Big) \\
    + \sum_{t=t0}^{8759}\sum_b\sum_{f} \phi^{CO_2,f} \cdot f_{f,t,b} \leq  \phi^{CO_2,el}_{t}\cdot\sum_{t=t0}^{8759} \bigg( \sum_b\sum_{est}\eta_{est}\\
    \cdot y_{t,est,b}^{exp} + \sum_b\sum_{g}y_{t,g,b}^{exp}\bigg) +Comp^{0\rightarrow t0}
\end{multline}
$Em^{0\rightarrow t0}$ and $Comp^{0\rightarrow t0}$ are the emission and compensation from the beginning of the year to the current timestep. 

This model is much longer to solve because of the number of timesteps in each iteration. In the other MPC, we chose $T_{MPC}=24$ timesteps from the beginning until the end while here it starts at 8760 and go down by one each time. Those implementation choices can be modulated depending on the computational load, by for example allowing to implement several hours instead of only the first one.

\section{Evaluation of the Operation Strategies} \label{op_res}

The different operation strategies presented in the previous sections are used to operate the systems resulting from the investment runs (and presented in section \ref{sec_inv_res}) in the years 2015, 2017 and 2018. We use a mipgap of 1\% and we use clusters for the perfect foresight and the receding horizon model in order to have reasonable solving time. We use 50 clusters for the perfect foresight (the same as in the investment runs) and 30 for the receding horizon. In the receding horizon case, we then have the $T_{mpc}$ hours from the current years and 30 clusters representing the remainder of the reference year instead of the whole reference year. Furthermore, for the receding horizon runs, we decide to implement the first 6 hours at each iteration instead of the first hour only in order to contain the computational time. In the MPC runs, we still only implement the first hour at each iteration.
For all runs, $T_{mpc}$ is set at 24 hours.

\begin{figure}[h]
    \centering
    \includegraphics[width=0.48\textwidth]{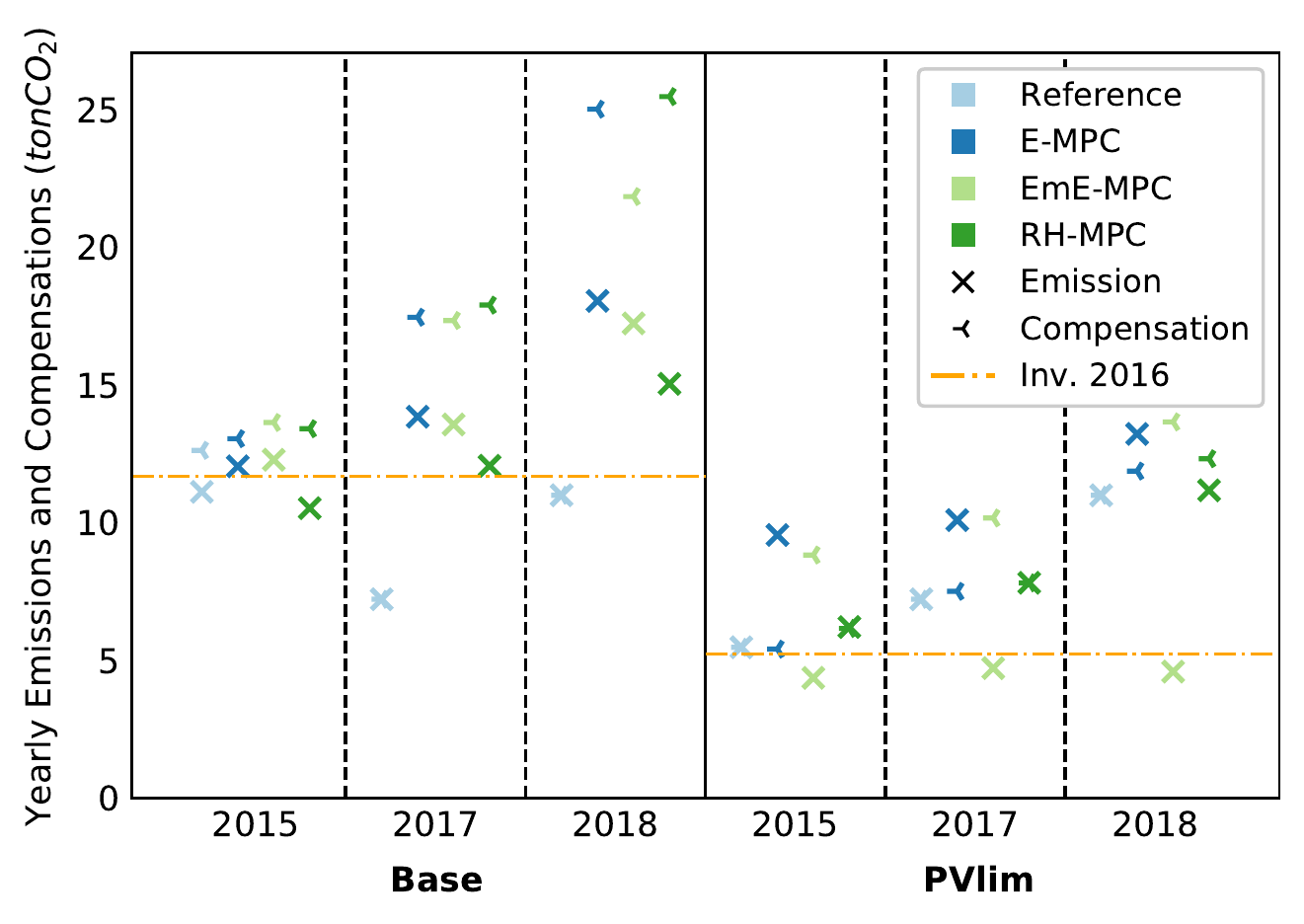}
    \caption{Emissions and Compensations from the Different Operation Strategies in the Different Years Considered}
    \label{op_em_comp}
\end{figure}

The yearly Emissions and Compensations are presented on Fig.\ref{op_em_comp}. The orange line, ``Inv. 2016" represents the level of emissions and compensations obtained in the investment run with year 2016. In the Base case, the emissions are always compensated. The energy system with its large amount of PV is quite passive and there is only the need to supply the heating load from the heat pumps and biomethane boiler. Even with the purely economic approach from E-MPC, the emission balance is satisfied. In years 2017 and 2018, the $CO_2$ factors (and to a lesser extent the spot prices and solar irradiances) are higher, making it harder for EmE-MPC to keep emissions at the level of the investment run.  Overall the RH-MPC approach gives the lowest emissions.

In the PVlim case, the system requires a more active management due to the lower amount of PV and the large CHP plant. 
It is not sufficient to operate the system in a purely cost optimal way because there is then not enough compensations, this is illustrated by the E-MPC approach. The EmE-MPC approach on the other hand keeps the emissions low and the compensations high. It manages to stay around the same level as in the investment run thanks to the penalization of deviating from the emissions and compensations resulting from the investment run. It manages to do so by using the CHP more even if it means dumping some of the heat produced. The RH-MPC approach gives again the best result. It manages to keep the total emissions and compensations close and they are always around the same level as in the Reference (perfect foresight approach). 

\begin{figure}[h]
    \centering
    \includegraphics[width=0.48\textwidth]{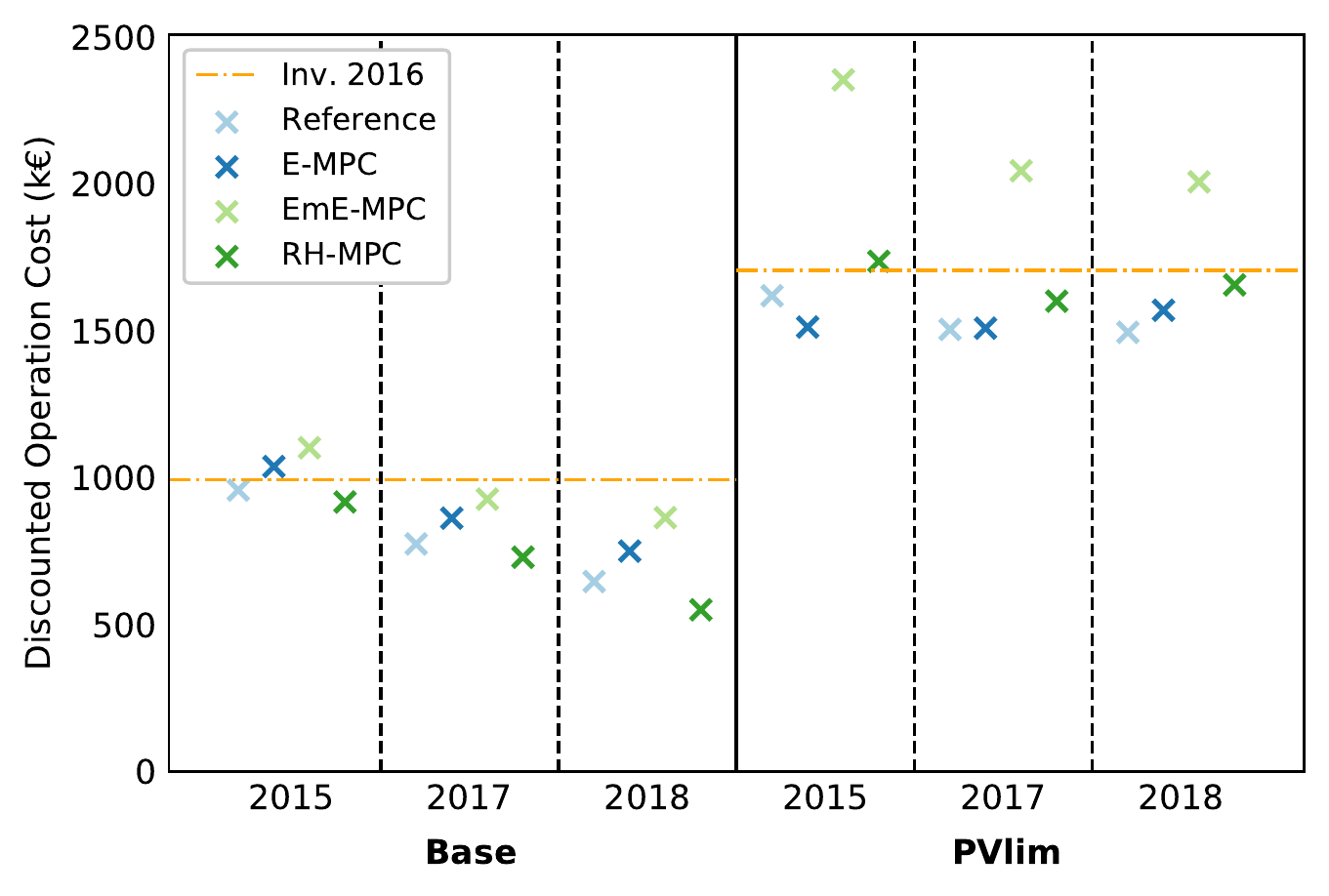}
    \caption{Discounted Operation Cost from the Different Operation Strategies in the Different Year Considered}
    \label{op_opc_comp}
\end{figure}

The total discounted operation costs are presented in Fig. \ref{op_opc_comp}. Note that the ``Inv. 2016" represents the operation costs from the investment runs and that the ``fictitious" penalization costs in the EmE-MPC are not included. 
In the Base case, the operation costs are lower for the years 2016 and 2017. This is partially due to the higher irradiance. The EmE-MPC has to follow the same pattern of emissions as in the investment run causing additional costs. RH-MPC has lower operation costs than the Reference, most likely because of the clustering. The Reference has 50 clusters for the year while the RH-MPC has actual data for 24 hours and clusters that are remade at each iteration giving a better representation of the year. 
In the PVlim case, the extra cost of maintaining the same emissions as in the investment run for EmE-MPC can be observed. They stem from the need to operate the costly CHP to reach the targets and avoid the penalization. The RH-MPC approach allows staying around the operation cost from the investment runs even though they are not as low as they could be (by comparison with the E-MPC for example).
%
%
%

\begin{figure}[h]
    \centering
    \includegraphics[width=0.48\textwidth]{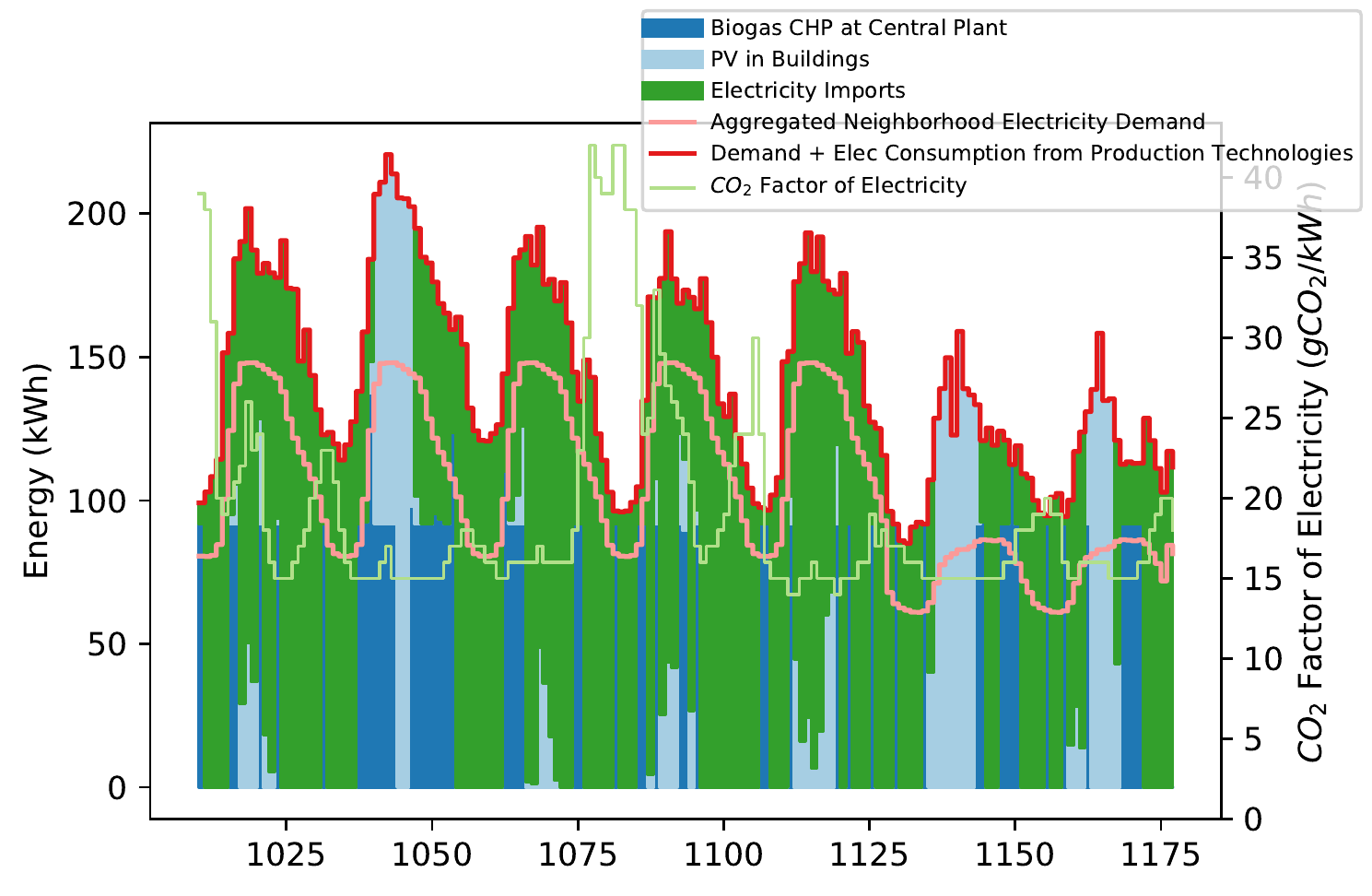}
    \caption{Origin of the Electricity Consumed in the Neighborhood and $CO_2$ Factor of Electricity in One Winter Week of 2018 in the PVlim E-MPC Case}
    \label{el_cons_E-MPC}
\end{figure}

\begin{figure}[h]
    \centering
    \includegraphics[width=0.48\textwidth]{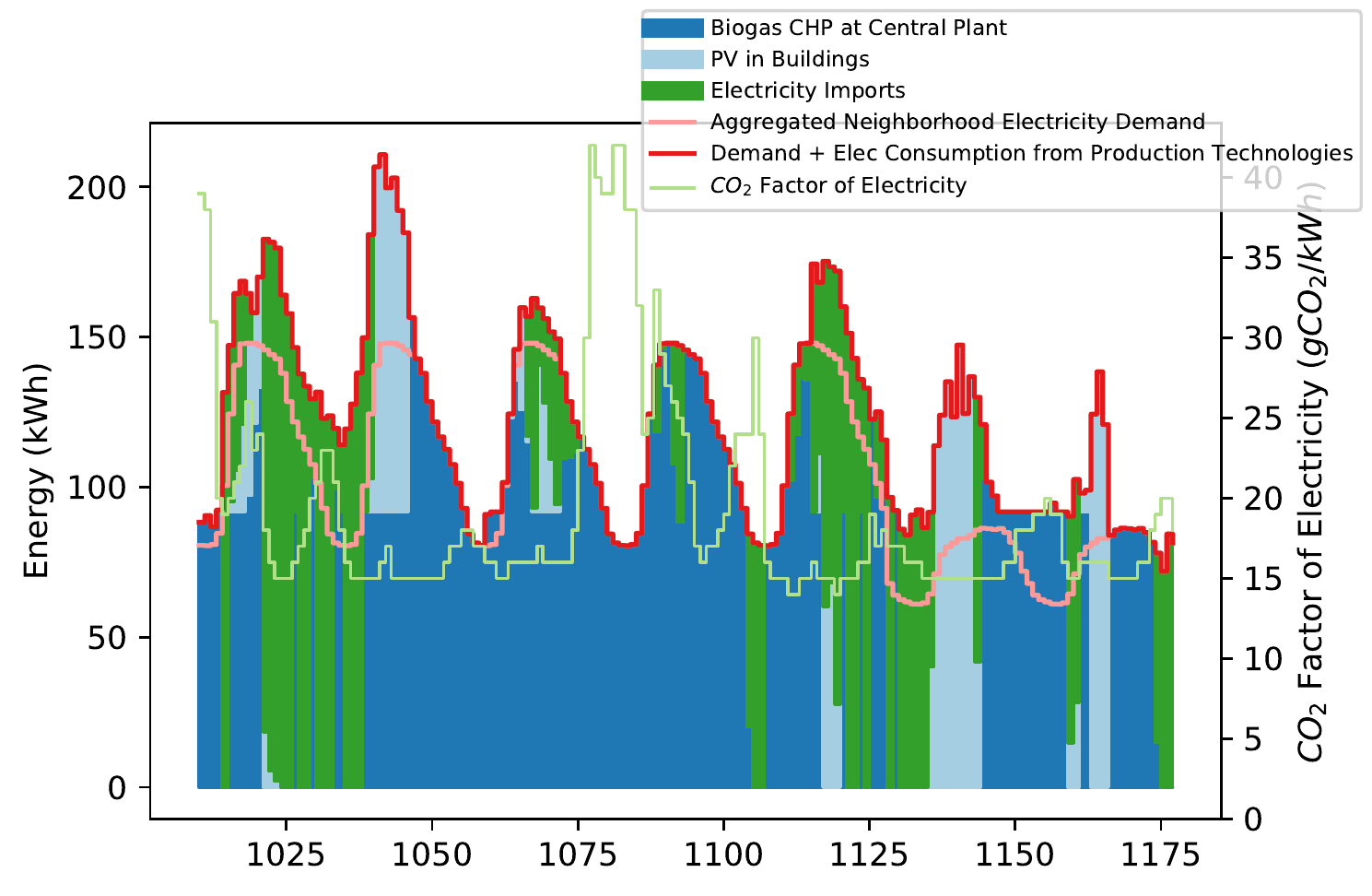}
    \caption{Origin of the Electricity Consumed in the Neighborhood and $CO_2$ Factor of Electricity in One Winter Week of 2018 in the PVlim EmE-MPC Case}
    \label{el_cons_EmE-MPC}
\end{figure}

\begin{figure}[h]
    \centering
    \includegraphics[width=0.48\textwidth]{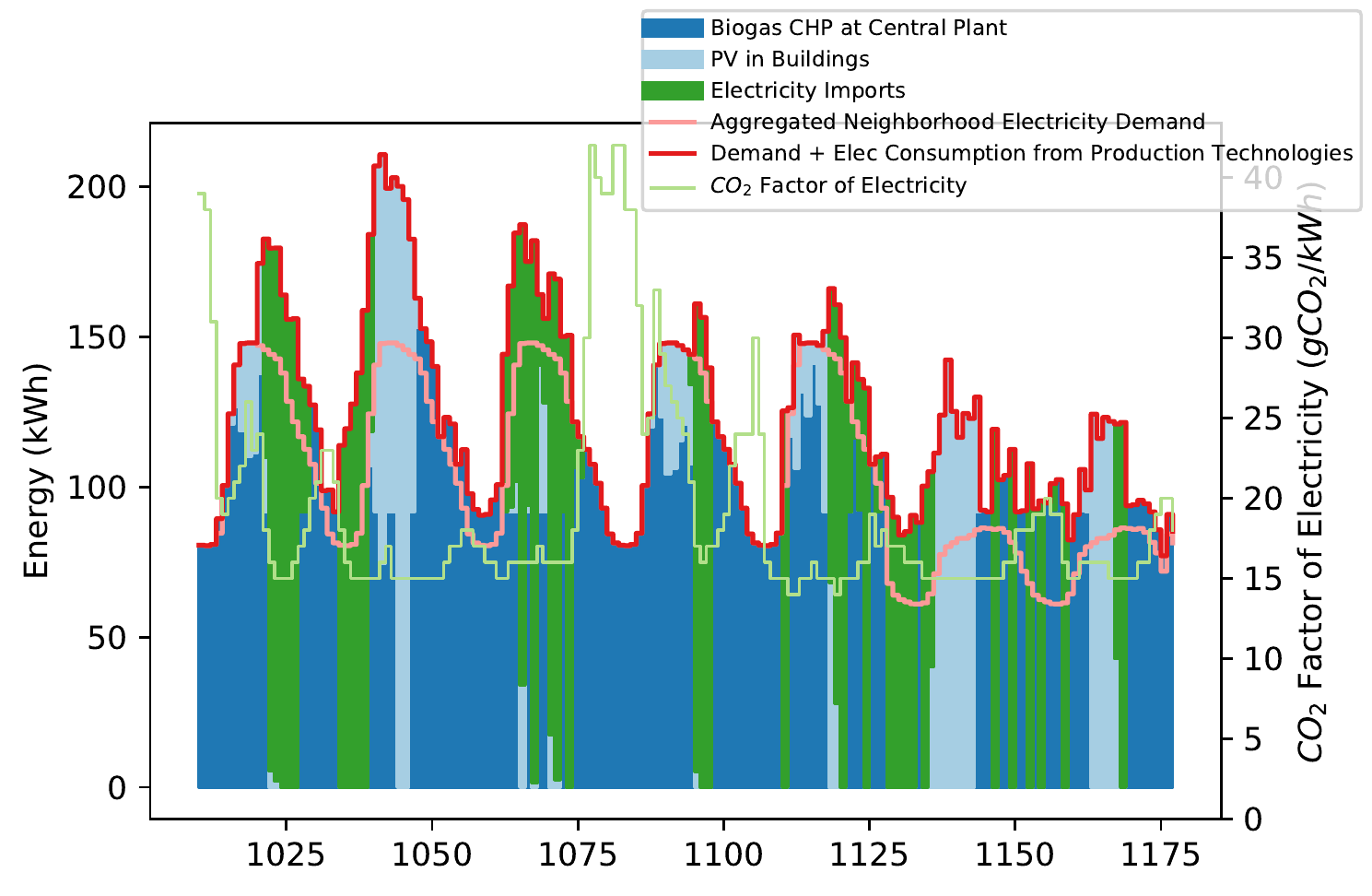}
    \caption{Origin of the Electricity Consumed in the Neighborhood $CO_2$ Factor of Electricity in One Winter Week of 2018 in the PVlim RH-MPC Case}
    \label{el_cons_RH-MPC}
\end{figure}

Fig. \ref{el_cons_E-MPC}, Fig. \ref{el_cons_EmE-MPC} and Fig. \ref{el_cons_RH-MPC} illustrate the differences in operation for one winter week in the year 2018 of the different operation strategies. This highlights the use of the CHP as a way to increase the compensation by exporting more and reduce emissions by importing less.
Fig. \ref{el_cons_RH-MPC} in particular amply illustrates the importance of the $CO_2$ factor of electricity in the choice of when to operate the CHP. The CHP is operated when the factor is high, i.e. when it is the most beneficial. In contrast, for the E-MPC, Fig. \ref{el_cons_E-MPC}, the operation is not so correlated to the $CO_2$ factor level. It is most likely more correlated to the spot price of electricity, which is in line with its purely economical approach.

%

In the EmE-MPC cases, especially in the PVlim, the CHP is operated to produce electricity that can be exported to the grid even though there is no need for the heat. This leads to a large amount of dumped heat. This is also the case to a lesser extent in the RH-MPC case. This could be solved by adding extra technological options dedicated to electricity production or CHP with heat to power ratios more in favor of electricity. In the Base case, there is also a great deal of electricity dumped due to the size of the grid connection. This is also the case in the investment runs. In the Base case in particular it is more cost efficient to over-invest in PV in order to be at maximum export for longer during the year, even if it means curtailing PV production at times.

\begin{figure}[h]
    \centering
    \includegraphics[width=0.48\textwidth]{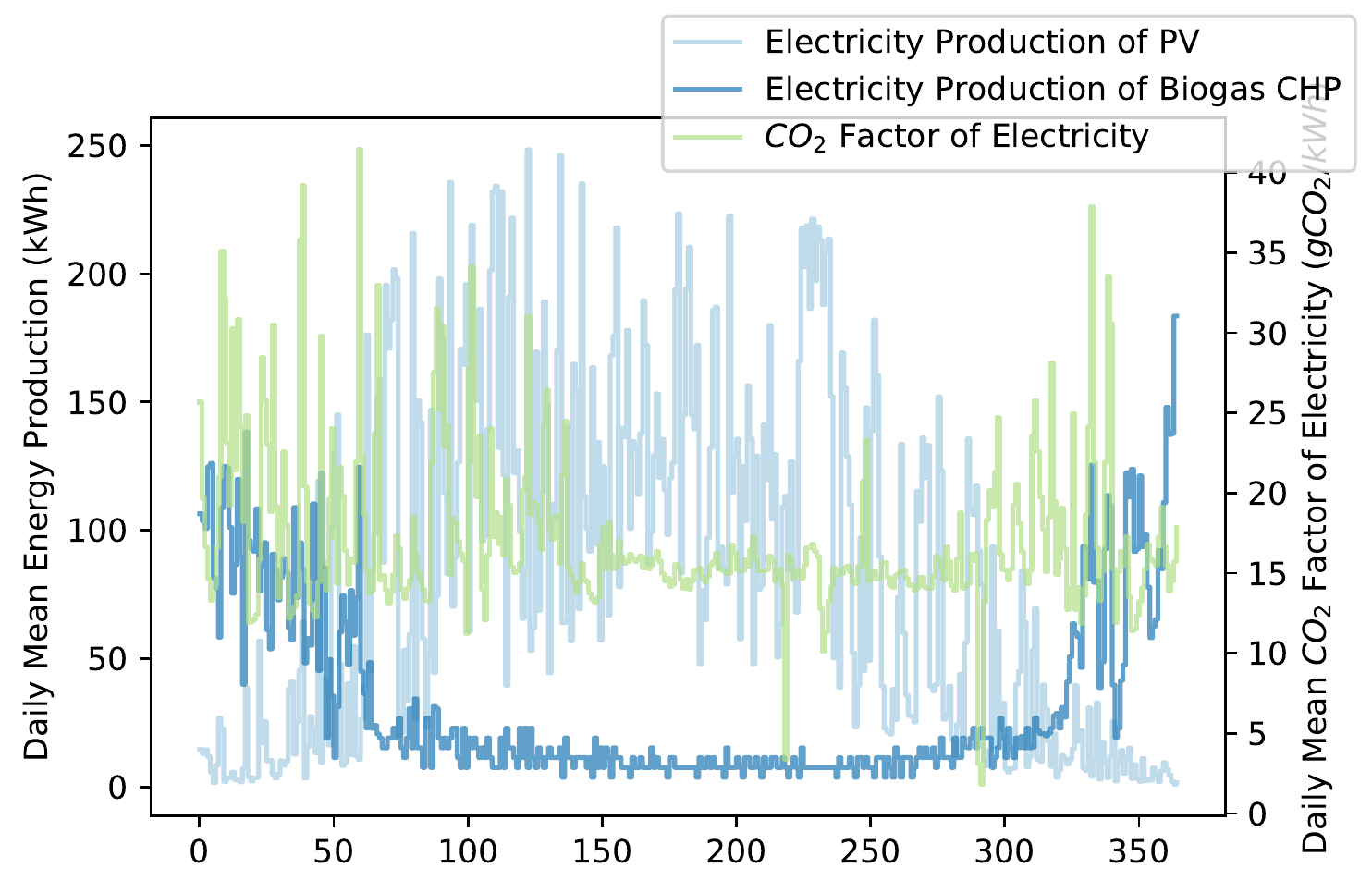}
    \caption{Daily Mean of $CO_2$ Factor of Electricity and Daily Mean Production of Electricity from the PV and the Biogas CHP in the PVlim 2015 RH-MPC Case}
    \label{prod_vs_factor}
\end{figure}

Both those electricity and heat dumps are linked to the $CO_2$ factor of electricity profiles. For example, Fig. \ref{prod_vs_factor} shows the daily mean production of electricity from the PV and from the Biogas CHP in the PVlim 2015 RH-MPC case. It highlights that the PV production, which cannot be controlled happens for a large part at times when the $CO_2$ factor is low. The CHP better matches the times of high factors due to its controllability in addition to the matching between the high winter thermal load and the high factors in the winter. This figures shows daily average for readability ease but note that the variations of $CO_2$ factors are larger and more frequent at the hourly level in the winter. In addition, year 2015 is the one with the lowest $CO_2$ factors of electricity levels and variability as can be seen in Fig. \ref{fig:box_co2}, Fig. \ref{fig:den_co2} and Fig. \ref{fig:dur_co2}.
One consequence of the RH-MPC method can be noted from this figure as well. In the last days of the year, the production of the CHP increases a lot. This is likely a result of the receding horizon approach. Due to the replacement of the reference year with the actual value for that year, there is a need to make up for the difference in emissions/compensations between what was expected and is possible with the actual data. 

\section{Limitations}

Several limitations should be kept in mind when it comes to the methodology and the interpretation of results. The first aspect to keep in mind is the effect of using clustering. The need to use clustering arises from the complexity of solving some of the models, in particular the investment model (Base and PVlim), the perfect foresight (PVlim) and the RH-MPC (PVlim) due to the binary variables. In order to keep the same conditions in all cases, clustering was used for all the appropriate cases (i.e. except the E-MPC and EmE-MPC that only consider a fixed ``short" horizon). This means the results are affected by the performance of the clustering and more information on this can be found in \cite{dimitri19}. Another parameter that was used for all cases for computational reason is a mipgap of 1\%. The PVlim investment in particular was converging very slowly below a mipgap of around 1\%.
Another limitation of the study is the years chosen. We chose years from 2015 to 2018 because these were the only ones where we could compute the hourly $CO_2$ factor of electricity from data available on the ENTSO-E transparency platform. This means that the interpretation of the results in a longer term setting is more uncertain. The profiles of the different timeseries could be at different levels or with different profiles in some decades and due to climate change.
EmE-MPC was only presented with one set of values for the parameters $\delta$, when in fact they would probably require fine tuning to be used in practice and have an effect that is just right and not be useless or too zealous. 
The use of clusters for RH-MPC makes it faster to solve but likely reduces its performance. Also this scheme requires more computation, even though the use of clustering partly alleviate this. This could be a problem in practice depending on the frequency of the optimization. 
A different approach to the RH-MPC would be to keep optimizing over a complete year through the iterations, without having it recede. This would remove the end of the year effect that was observed and give an homogeneous solving time throughout the year.

\section{Conclusion}

In the concept of ZEN there is a need for a better transition between the design recommendation from investment tools that assume a certain operation and the way the energy system would actually be operated. In particular the strong requirement on emissions cannot be considered in the same way in the operation and in the investment process. In this paper we suggested and compared different operation approaches and their performance in terms of operation cost and emissions/compensations. The investment tool ZENIT was first used to create designs of ZEN energy systems in cases where the amount is and is not limited. We then compared the performance of four approaches in operating those systems in different years. The first one used as a reference assumes perfect foresight of the year and is used as a reference; the E-MPC approach represents a purely economical operation of the neighborhood; the EmE-MPC approach expends the E-MPC by including a penalization of deviating from emission and compensation targets and the RH-MPC approach uses a receding horizon and a complete year as a way to maintain the annual zero emission constraint in the short-term operation optimization.
We also look into the variations between data from different years and how this affects the actual costs. Indeed in the investment run we use a reference year and expect the operation cost, emissions and compensation for the actual operation of the ZEN to even themselves out between years.
The results show that with a system strongly based on PV, the zero emission requirement can be met without any additional specific operation method. However in systems including technologies using carbon-intensive sources or systems where one of the source is expensive to operate (such as the CHP in the PVlim case) the need for a more active operation and for accounting emissions and compensations in some way is greater. To this end, the proposed RH-MPC appears to be the most promising operation strategy. The EmE-MPC method perfomred less well, but better tuning of penalization cost parameters could make also make this a viable solution.
This study could be expanded in the future by considering other approaches for the operation of the neighborhood and also by considering a ZEN energy system that includes more carbon-intensive sources, for instance by having a lower requirement for the compensation and only partially compensating emissions partly.

\section*{Acknowledgment}
This article has been written within the Research Center on Zero Emission Neighborhoods in Smart Cities (FME ZEN). The author gratefully acknowledges the support from the ZEN partners and the Research Council of Norway.

\bibliographystyle{IEEEtran}
\bibliography{biblio}


%

\onecolumn
\newpage
\appendices

\section{Technology Data}\label{annex:data}
The data for those technologies come from the Danish Energy Agency and Energinet\footnote{
\url{https://ens.dk/en/our-services/projections-and-models/technology-data}}.
\begin{threeparttable}[h]
    \centering
    \caption{Data Of Technologies Producing Heat and/or Electricity}
    \scriptsize 	
    \begin{tabular}{l r r r r r r r r r r r}
    \hline
        Tech. & $\eta_{th}$ & Fix. Inv. Cost  & Var. Inv. Cost  & $\alpha_i$ & Min. Cap. & Annual O\&M Costs & Lifetime & Fuel & $\alpha_{CHP}$ & El. & Heat\\
         & (\%) & (\euro) &  (\euro/kW) & (\% Inst. Cap.) & (kW)  &  (\% of Var Inv. Cost) & (year) & &  \\
        \hline
        \multicolumn{4}{l}{\textbf{At building level}} & & & & &  & & \\
        \hline
        PV\tnote{1} & & 0 & 730 & 0 & 50 & 1.42 & 35 &  & & 1 & 0\\
        ST\tnote{2} & 70 & 28350 & 376 & 0 & 100 & 0.74 & 25 &  & & 0 & 1 \\
        ASHP\tnote{3} & f($T_t$) & 42300 & 247 & 0 & 100 & 0.95 & 20 & Elec. &  & 0 & 1 \\
        GSHP\tnote{4} & f($T_t$) & 99600 & 373 & 0 & 100 & 0.63 & 20 & Elec. &  & 0 & 1 \\
        Boiler\tnote{5} & 85 & 32200 & 176 & 30 & 100 & 2.22 & 20 & Wood Pellets  & &0 & 1 \\
        Heater & 100 & 15450 & 451 & 0 & 100 & 1.18 & 30 & Elec.  & &0 & 1 \\
        Boiler & 100 & 3936 & 52 &20 & 35 & 2.99 & 25& Biomethane & & 0 &1 \\
        \hline
        \multicolumn{4}{l}{\textbf{At neighborhood level}} & & & &  & & & \\
        \hline
        CHP\tnote{6} & 47 & 0 & 1035 & 50 & 200 & 1.03 & 25 & Biogas &  1.09 & 1 & 1\\
        CHP & 98 & 0 & 894 & 20 & 1000 & 4.4 & 25 & Wood Chips & 7.27 & 1 & 1 \\
        CHP & 83 & 0 & 1076 & 20 & 1000& 4.45 & 25 & Wood Pellets & 5.76 & 1 & 1 \\
        Boiler\tnote{7} & 115 & 0 & 680 & 20 & 1000 & 4.74 & 25 & Wood Chips   & & 0 & 1 \\
        Boiler\tnote{7} & 100 & 0 & 720 & 40 & 1000& 4.58 & 25 & Wood Pellets   & & 0 & 1 \\
        CHP\tnote{8} & 66 & 0 & 1267 & 10 & 10 & 0.84 & 15 & Wood Chips &  3 & 1 & 1 \\
        Boiler\tnote{9} & 58 & 0 & 3300 & 70 & 50 & 5 & 20 & Biogas &  &0 & 1\\
        GSHP\tnote{4} & f($T_t$) & 0 & 660 & 010 & 1000& 0.3 & 25 & Elec. &  & 0 & 1\\
        Boiler & 99 & 0 & 150 & 5 & 60 & 0.71 & 20 & Elec. &   & 0 & 1 \\
        Boiler & 100 & 0& 60 & 15 & 500 & 3.25 & 25 & Biogas &  & 0&1\\
    
    \hline    
    \end{tabular}
    \begin{tablenotes}
    \scriptsize
    \setlength{\columnsep}{0.8cm}
    \setlength{\multicolsep}{0cm}
    \begin{multicols}{3}
     \item[1] Area Coefficient: 5.3 $m^2/kW$
     \item[2] Area Coefficient: 1.43 $m^2/kW$
     \item[3] Air Source Heat Pump
     \item[4] Ground Source Heat Pump
     \item[5] Automatic stoking of pellets
     \item[6] Gas Engine
     \item[7] HOP
     \item[8] Gasified Biomass Stirling Engine Plant
     \item[9] Solid Oxyde Fuel Cell (SOFC)
    
    \end{multicols}
    \end{tablenotes}
    \label{tab:techdata}
\end{threeparttable}
\vspace{0.3cm}

The data for prices of fuels come from different sources. For the wood pellets and wood chips, they come from the Norwegian Bioenergy Association\footnote{\url{http://nobio.no/wp-content/uploads/2018/01/Veien-til-biovarme.pdf}}. The data for the biogas and biomethane come from the European Biogas Association\footnote{ \url{https://www.europeanbiogas.eu/wp-content/uploads/2019/07/Biomethane-in-transport.pdf}}.

The data for $CO_2$ factor of fuels come from a report from Cundall\footnote{\myurl}.

\begin{table}[h]
    \centering
    \caption{Data of Fuels}
    \label{tab:fueldata}
    \begin{tabular}{l r r}
    \hline
        Fuel & Fuel Cost ($\euro/kWh$) & $CO_2$ factor ($gCO_2/kWh$)\\
         \hline
        Electricity & $f(t)$ & $f(t)$ \\
        Wood Pellets & 0.03664 & 40 \\
        Wood Chips & 0.02592 & 20 \\
        Biogas & 0.07 & 0 \\
        Biomethane & 0.07 & 100\\
         \hline
    \end{tabular}
\end{table}
\begin{threeparttable}[h]
    \centering
    \caption{Data of Storage}
    \scriptsize
    \begin{tabular}{l r r r r r r }
    \hline
        Index & One way eff. & Inv. Cost & O\&M Cost & Lifetime & Min. Cap. & Charge / Discharge rate \\
         & (\%) & (\euro/kWh) & (\% of Inv. Cost) & (year) & (kWh) & (\% of Cap) \\
         \hline
         \multicolumn{7}{l}{Battery} \\
         \hline
         1\tnote{1} & 95 & 577 & 0 & 10 & 13.5 & 37 \\
         2\tnote{2} & 938 & 500 & 0 & 15 & 210 & 23 \\
         3\tnote{3} & 95 & 432 & 0 & 20 & 1000 & 50 \\
         \hline
         \multicolumn{7}{l}{Heat Storage}\\
         \hline
         1\tnote{4} & 95 & 75 & 0 & 20 & 0 & 20 \\
         2\tnote{3} & 98 & 3 & 0.29 & 40 & 45 000 & 1.7 \\
         \hline
    \end{tabular}
    \begin{tablenotes}
    \scriptsize
    \setlength{\columnsep}{0.8cm}
    \setlength{\multicolsep}{0cm}
    \begin{multicols}{2}
     \item[1] Based on Tesla Powerwall
     \item[2] Based on Tesla Powerpack
     \item[3] Based on Danish energy agency data
     \item[4] Same data are used for the heat storage at the building or neighborhood level and for both SH and DHW
    \end{multicols}
    \end{tablenotes}
    \label{tab:stordata}
\end{threeparttable}

%

\ifCLASSOPTIONcaptionsoff
  \newpage
\fi

\end{document}